\documentclass[preprint]{imsart}
\usepackage[margin=1in]{geometry}

\usepackage{cite}

\usepackage{graphicx}
\DeclareGraphicsExtensions{.pdf}
\graphicspath{{./figures/}}

\usepackage{tabularx}
\usepackage{array}
\usepackage{dcolumn}
\newcolumntype{d}[1]{D{.}{.}{#1}}
\newcolumntype{.}{D{.}{.}{-1}}
\usepackage{booktabs}
\usepackage{threeparttable}
\usepackage{rotating}
\usepackage{multirow}

\usepackage{textcomp}

\usepackage[caption=false,font=footnotesize]{subfig}

\usepackage[hyphens]{url}

\usepackage{fixltx2e}

\usepackage[cmex10]{amsmath}
\usepackage{mathtools}
\usepackage{bbm}

\newcommand{\argmin}[1]{\underset{#1}{\operatorname{argmin}}\ }
\newcommand{\funct}[1]{\mathrm{#1}}
\newcommand{\transpose}[1]{#1^\mathrm{T}}

\newcommand{\norm}[2]{\left\| #1 \right\|_{#2}}

\usepackage{amssymb}
\newcommand{\Real}{{\mathbb R}}
\newcommand{\Complex}{{\mathbb C}}

\usepackage{algpseudocode}
\usepackage{eqparbox}
\algrenewcommand{\algorithmiccomment}[1]{\hfill // \eqparbox{COMMENT}{#1}}

\newsavebox{\ieeealgbox}
\newenvironment{boxedalgorithmic}
  {\begin{lrbox}{\ieeealgbox}
   \begin{minipage}{\dimexpr\textwidth-2\fboxsep-2\fboxrule}
   \begin{algorithmic}}
  {\end{algorithmic}
   \end{minipage}
   \end{lrbox}\noindent\fbox{\usebox{\ieeealgbox}}}

\begin{document}

\begin{frontmatter}
\title{L1 Splines for Robust, Simple, and Fast Smoothing of Grid Data}

\begin{aug}
  \author{
    \fnms{Mariano} \snm{Tepper}
    \ead[label=e1]{mariano.tepper@duke.edu}
  }

  \author{
    \fnms{Guillermo} \snm{Sapiro}
    \ead[label=e2]{guillermo.sapiro@duke.edu}
  }

  \address{
    Department of Electrical and Computer Engineering, Duke University.\footnote{This work was partially done while the authors were with the Department of Electrical and Computer Engineering, University of Minnesota.}\\
    \printead{e1,e2}
  }
  
  \runauthor{Tepper and Sapiro}
  \runtitle{L1 Splines for Robust, Simple, and Fast Smoothing of Grid Data}

\end{aug}

\begin{abstract}
Splines are a popular and attractive way of smoothing noisy data. Computing splines involves minimizing a functional which is a linear combination of a fitting term and a regularization term. The former is classically computed using a (weighted) L2 norm while the latter ensures smoothness. Thus, when dealing with grid data, the optimization can be solved very efficiently using the DCT. In this work we propose to replace the L2 norm in the fitting term with an L1 norm, leading to automatic robustness to outliers. To solve the resulting minimization problem we propose an extremely simple and efficient numerical scheme based on split-Bregman iteration combined with DCT. Experimental validation shows the high-quality results obtained in short processing times.
\end{abstract}

\begin{keyword}
\kwd{Splines}
\kwd{L1 fitting}
\kwd{split-Bregman}
\kwd{grid data}
\end{keyword}

\end{frontmatter}

\maketitle

\section{Introduction}

Smoothing a dataset consists in finding an approximating function that captures important patterns in the data, while disregarding noise or other fine-scale structures. Let $y \in \Real^{n_1 \times \dots \times n_m} \rightarrow \Complex$ be an $m$-dimensional discrete signal, where $n_j$ ($1 \leq j \leq d$) is the domain of $y$ along the $j$-th dimension. We can model $y$ by
\begin{equation}
    y = \hat{y} + r ,
    \label{eq:model}
\end{equation}
where $r$ represents some noise and $\hat{y}$ is a smooth function. A very common regularization choice is to enforce $C^2$ continuity, in which case $\hat{y}$ is called a (cubic) spline. Smoothing $y$ relies upon finding the best estimate of $\hat{y}$ under the proper smoothness and noise assumptions.
We can approximate $\hat{y}$ by minimizing an objective functional
\begin{equation}
    \funct{F}(z) = \funct{R}_y (z) + s \funct{P}(z) ,
\end{equation}
where $\funct{R}_y (z)$ is a data fitting term, defined by the distribution of $r$, $\funct{P}(z)$ is a regularization term, and $s$ is a scalar that determines the balance between both terms. Such scalar parameter can be automatically derived using Bayesian or MDL techniques, as will be later shown.

For clarity, we will describe in depth the case $m=1$, where we have $n=n_1$ samples. We extend the results later to the general $m$-dimensional case.

Let us begin by explaining the smoothing term.
The $C^2$ continuity requirement leads to define $\funct{P}(z) = \norm{Dz}{2}^ 2$, $D$ being a discrete \emph{second-order} differential operator, defined $\forall i, 2 \leq i \leq n-1$, by
\begin{align}
    D_{i, i-1} &= \tfrac{2}{h_{i-1} (h_{i-1} + h_{i})} ,\\
    D_{i, i} &= \tfrac{-2}{h_{i-1} \ h_{i}} ,\\
    D_{i-1, i} &= \tfrac{2}{h_{i} (h_{i-1} + h_{i})} ,
\end{align}
where $h_i$ represents the step, or sampling rate, between $y_i$ and $y_{i+1}$.
Assuming repeating border elements, that is, $y_0 = y_1$ and $y_{n+1}  =y_n$, gives $D_{1, 1} = -D_{1, 2} = -1/h_1^2$, and $D_{n, n-1} = -D_{n, n} = -1/h_{n-1}^2$.

Regarding the fitting term, the classical assumption is that the noise $r$ in Equation~(\refeq{eq:model}) has Gaussian distribution with zero mean and unknown variance, which leads to setting $\funct{R}_y (z) = \norm{z - y}{2}^2$. 
Smoothing then can be formulated as the least-squares regression
\begin{equation*}
    \hat{y} = \argmin{z} \norm{z - y}{2}^2 + s \norm{Dz}{2}^ 2 .
\end{equation*}
For clarity, we call the estimate $\hat{y}$ obtained with this method an L2 spline.

It is a well known fact that least squares estimates for regression models are highly non-robust to outliers. Although there is no agreement on a universal and formal definition of an outlier, it is usually regarded as an observation that does not follow the patterns in the data. Notice that smoothing should produce an estimate $\hat{y}$ taking into account only important patterns, that is, the inliers, in the data. In this sense, the L2 formulation cannot correctly handle outliers by itself.

In order to solve this problem, we propose to take a different assumption on the distribution of the noise $r$ in Equation~(\refeq{eq:model}). By choosing a distribution with fatter tails than the Gaussian distribution, the derived estimator will correctly handle outliers. We thus assume that $r$ follows a Laplace distribution with zero mean and unknown scale parameter, a common practice in other problems as we will further discuss below. This leads to the fitting term $\funct{R}_ y(z) = \norm{z - y}{1}$ and the regression then becomes
\begin{equation*}
    \hat{y} = \argmin{z} \ \norm{z - y}{1} + s \norm{Dz}{2}^2 .
\end{equation*}
We call the estimate $\hat{y}$ obtained with this formulation an L1 spline.

Let us point out that the use of L1 fitting terms for solving inverse problems is not new.
For example, in 2001, Nikolova proved the theoretic pertinence of using L1 fitting terms for image denoising~\cite{nikolova02}. 
Other interesting works have addressed this approach for total variation image denoising~\cite{alliney97,chan2004,aujol06,nikolova12} or total variation optical flow~\cite{zach07,wedel09,raket11} (this robustness-to-outliers type of ideas was previously introduced in the context of optical flow by Black and Anandan~\cite{black91}). Recall that the total-variation regularization term involves first-order derivatives, while the proposed L1 splines, on the other hand, involve second-order derivatives.

Regularly sampled signals are extremely common in practice, and their analysis becomes easier and faster. In particular, we follow the most common choice when dealing with discrete $m$-dimensional data, which is assuming a ``rectangular'' Cartesian sampling pattern. When the sampling is isotropic, i.e., ``square,'' we refer to this type of data as grid data.

We developed an iterative algorithm for computing L1 splines, based on split-Bregman iteration~\cite{goldstein09}, that is specially suited for the case of grid data. This algorithm is extremely fast, both in running time and in the number of iterations until convergence. It is also outstandingly simple, making the implementation completely straightforward.

The remainder of the paper is structured as follows. In Section~\ref{sec:L2splines}, we overview a fast algorithm to compute L2 splines and robust L2 splines, a modification of the least squares regression that allows to handle outliers. In Section~\ref{sec:L1splines}, we present the proposed algorithm for computing L1 splines. Then, in Section~\ref{sec:results}, we show results obtained with L1 splines, systematically outperforming its L2 and robust L2 counterparts in the presence of outliers. We also show that the proposed computational algorithm is very efficient. Finally, in Section~\ref{sec:conclusions} we provide some concluding remarks.

\section{Smoothing splines}
\label{sec:L2splines}

As aforementioned, the classical assumption is that the noise $r$ in Equation~(\refeq{eq:model}) follows a Gaussian distribution with zero mean and unknown variance. This leads to solve the least-squares regression
\begin{equation}
    \hat{y} = \argmin{z} \norm{z - y}{2}^2 + s \norm{Dz}{2}^ 2 .
    \label{eq:LSformulation}
\end{equation}
Since both terms are differentiable, we obtain
\begin{equation}
     \hat{y} = (I + s \transpose{D} D)^{-1} y .
     \label{eq:LSsolution}
\end{equation}

Garcia proposed a very efficient method for dealing with regularly sampled data~\cite{garcia10}. Assuming that the data are equally spaced, that is, without loss of generality $\forall i, h_i = 1$, we obtain
\begin{equation}
D = 
\begin{pmatrix}
-1 & 1 \\
1 & -2 & 1 \\
& \ddots & \ddots & \ddots \\
&& 1 & -2 & 1 \\
&&& 1 & -1
\end{pmatrix} .
\end{equation}
An eigendecomposition of $D$ yields $D = U \varLambda \transpose{U}$, where $\varLambda$ is a diagonal matrix containing the eigenvalues of $D$, given by~\cite{yueh05}
\begin{equation}
    \varLambda_{i,j} =
    \begin{cases}
        -2 + 2 \cos((i-1) \pi / n) , & \text{if $i=j$;} \\
        0 , & \text{otherwise.}
    \end{cases}
\end{equation}
Since $U$ is a unitary matrix, we can write Equation~(\refeq{eq:LSsolution}) as
\begin{equation}
     \hat{y} = U (I + s \varLambda^2)^{-1} \transpose{U} y .
     \label{eq:LSsolutionDecomposed}
\end{equation}
Let us define the matrix $\varGamma = (I + s \varLambda^2)^{-1}$. Trivially,
\begin{equation}
    \varGamma_{i, j} =
    \begin{cases}
     \left[ 1 + s (-2 + 2 \cos((i-1) \pi / n))^2 \right]^{-1} , & \text{if $i=j$;} \\
     0 , & \text{otherwise.}
     \end{cases}
     \label{eq:varGamma}
\end{equation}
Following Strang~\cite{strang99} and Garcia~\cite{garcia10}, let us observe that $\transpose{U}$ is a DCT-II matrix and $U$ is an inverse DCT-II matrix. Then, Equation~(\refeq{eq:LSsolutionDecomposed}) can be expressed as
\begin{equation}
    \hat{y} = \funct{DCT^{-1}} (\varGamma \ \funct{DCT} (y)) ,
    \label{eq:LSsolutionDCT}
\end{equation}
where $\funct{DCT}(\cdot)$ and $\funct{DCT^{-1}}(\cdot)$ stand for the DCT-II and inverse DCT-II functions. Equation~(\refeq{eq:LSsolutionDCT}) provides a fast and simple algorithm for computing L2 splines.

\subsection{Robust estimation.}
\label{sec:robustL2}

Often in practice there are in $y$ some values $y_i$ that could not be observed (or recorded) for some reason.
We would like to be able to handle such cases in such a way that the missing values are inferred from the ones that can be observed.
Let $W$ be an $n \times n$ diagonal matrix such that $W_{i, i}$ represents a weight assigned to observation $i$. $W$ is defined by
\begin{equation}
    W_{i, i} = 
    \begin{cases}
        0 &\text{if datapoint $i$ is missing;} \\
        \rho &\text{otherwise.}
    \end{cases}
\end{equation}
where $\rho$ is some arbitrary constant in $(0, 1]$; in practice, and without loss of generality, we set $\rho=1$. We can then solve
\begin{equation}
    \hat{y} = \argmin{z} \norm{ W^{1/2} (z - y)}{2}^2 + s \norm{Dz}{2}^ 2 ,
    \label{eq:WLSformulation}
\end{equation}
which will simply omit the missing points from the computation of the residual while the regularizer will still have a smoothing effect over both present and missing points. Equation~(\refeq{eq:WLSformulation}) acts as an impainting algorithm, filling the missing values in such a way that continuity between filled values and smoothed ones is preserved.
The minimization of Equation~(\refeq{eq:WLSformulation}) gives
\begin{equation}
    (I + s \transpose{D}D) \hat{y} = W (y - \hat{y}) + \hat{y} .
\end{equation}
This leads to the iterative procedure
\begin{equation}
    \hat{y}^{k+1} = (I + s \transpose{D}D)^{-1} \left( W \left( y - \hat{y}^k  \right) + \hat{y}^k \right) ,
\end{equation}
which, similarly to Equation~(\refeq{eq:LSsolutionDCT}), becomes
\begin{equation}
    \hat{y}^{k+1} = \funct{DCT}^{-1} \left( \varGamma \ \funct{DCT} \left( W \left(y - \hat{y}^k \right) + \hat{y}^k \right) \right) .
    \label{eq:WLSsolutionDCT}
\end{equation}

On a different note, real data often present observations that lie abnormally far from their ``true'' value, i.e., that do not appear to follow the pattern of the other data points.
The main drawback of the penalized least squares formulation Equation~(\refeq{eq:LSformulation}) is its sensitivity to these outliers.
To address this issue, weights can be assigned to every point, as in Equation~(\refeq{eq:WLSformulation}), such that outliers exert less influence during the estimation process.
In this case, the weights are iteratively refined during the estimation process using robust estimators for the mean and variance of the data. Defining these estimators is a complex problem by itself. For details about how $W$ can be set and updated for added robustness to outliers, refer to Garcia's work~\cite{garcia10}.

\section{L1 splines}
\label{sec:L1splines}

In this section we introduce a different splines formulation in order to handle outliers in the data. We assume that the noise $r$ in Equation~(\refeq{eq:model}) follows a Laplace distribution with zero mean and unknown scale parameter, which leads to $\funct{R}_ y(z) = \norm{z - y}{1}$. The regression then becomes
\begin{equation}
    \min_{z} \ \norm{z - y}{1} + s \norm{Dz}{2}^2 .
    \label{eq:L1formulation}
\end{equation}

Goldstein and Osher~\cite{goldstein09} proposed a very elegant and efficient algorithm for solving the L1 constrained problem
(related to a number of very efficient optimization algorithms, e.g., see~\cite{combettes11})
\begin{equation}
    \min_u \norm{\Phi(u)}{1} + \funct{H}(u) .
\end{equation}
For this, they consider the equivalent problem
\begin{equation}
    \min_u \norm{d}{1} + \funct{H}(u) \quad \text{s.t.} \quad d = \Phi(u) ,
\end{equation}
which they first convert it into the unconstrained problem
\begin{equation}
    \min_u \norm{d}{1} + \funct{H}(u) + \tfrac{\lambda}{2} \norm{d - \Phi(u)}{2}^2 .
\end{equation}
In this form, the penalty function does not accurately enforce the constraint for small $\lambda$. The constraint is enforced by letting $\lambda \rightarrow \infty$.
However, another solution for this new formulation is found by using the following two-phase algorithm
\begin{align}
    \left( u^{k+1}, d^{k+1} \right) &= \argmin{u} \norm{d}{1} + \funct{H}(u) +
    \tfrac{\lambda}{2} \norm{d - \Phi(u) - b^k}{2}^2 ,\\
    b^{k+1} &= b^k + \left( \Phi(u^{k+1}) - d^k \right) .
\end{align}
This algorithm is often denoted in the literature as split-Bregman iteration. This class of algorithms has several nice theoretical properties and has successfully been applied to several problems in practice such as image restoration~\cite{osher05}, image denoising~\cite{xu07}, compressed sensing~\cite{yin08}, and image segmentation~\cite{goldstein10}; see also~\cite{combettes11} and references therein.

We use this technique for solving Equation~(\refeq{eq:L1formulation}). We begin by setting $\Phi(z) = z-y$ and $H(z) = s \norm{Dz}{2}^2$, which leads to the problem
\begin{equation}
    \min_{z, d} \ \norm{d}{1} + s \norm{Dz}{2}^ 2  \quad \text{s.t.} \quad d = z - y .
\end{equation}
We then transform it into the unconstrained form
\begin{equation}
    \min_{z, d} \ \norm{d}{1} + s \norm{Dz}{2}^ 2 + \tfrac{\lambda}{2} \norm{d - z + y}{2}^2 ,
    \label{eq:L1formulationExtended}
\end{equation}
and the Bregman iteration simply takes the form
\begin{align}
    \left( z^{k+1}, d^{k+1} \right) & = \argmin{z, d} \ \norm{d}{1} + s \norm{Dz}{2}^ 2 +
    \tfrac{\lambda}{2} \norm{d - z + y - b^k}{2}^2 ,
    \label{eq:splitBregman1}\\
    b^{k+1} & = b^k + (z^{k+1} - y - d^{k+1}) . \label{eq:splitBregman2}
\end{align}

Because of the splitting of the L1 and L2 components in the functional~(\refeq{eq:splitBregman1}), we can perform this minimization efficiently by iteratively minimizing with respect to $z$ and $d$ separately,
\begin{align}
    z^{k+1} &= \argmin{z} \ s \norm{Dz}{2}^ 2 + \tfrac{\lambda}{2} \norm{d^k - z + y - b^k}{2}^2 \label{eq:splitBregman1_1} , \\
    d^{k+1} &= \argmin{d} \ \norm{d}{1} + \tfrac{\lambda}{2} \norm{d - z^{k+1} + y - b^k}{2}^2 \label{eq:splitBregman1_2} .
\end{align}
For minimizing Equation~(\refeq{eq:splitBregman1_1}) we set $\tilde{y} = d^k + y - b^k$ and $\tilde{s} = 2s / \lambda$. We obtain
\begin{equation}
    \min_{z} \ \tilde{s} \norm{Dz}{2}^ 2 + \norm{\tilde{y} - z}{2}^2 .
\end{equation}
This is a classical L2 spline and can be minimized using Equation~(\refeq{eq:LSsolutionDCT}), as already explained.
The optimal value of $d$ in Equation~(\refeq{eq:splitBregman1_2}) can be explicitly computed using shrinkage operators,
\begin{align}
    d^{k+1} &= \argmin{d} \ \norm{d}{1} + \tfrac{\lambda}{2} \norm{d - z^{k+1} + y - b^k}{2}^2 \nonumber \\
    &= \funct{Shrink} (z^{k+1} - y + b^k , 1/\lambda) ,
\end{align}
where
\begin{align}
    \funct{Shrink}(v, \gamma) = 
    \begin{pmatrix}
    \funct{shrink}(v_1, \gamma) \\
    \vdots \\
    \funct{shrink}(v_j, \gamma) \\
    \vdots \\
    \funct{shrink}(v_m, \gamma)
    \end{pmatrix}
    \intertext{and}
    \funct{shrink}(x, \gamma) = \frac{x}{|x|} \max(|x|-\gamma, 0) .
    \label{eq:shrinkage}
\end{align}
We thus obtain a very efficient algorithm for computing L1 splines, combining DCT and shrinkage operators.

On a different note let us mention that Equation~(\refeq{eq:L1formulationExtended}) can also be interpreted~\cite{mateos12} as a relaxation of
\begin{equation}
    \min_{z, d} \norm{d}{0} + s \norm{Dz}{2}^2 + \tfrac{\lambda}{2} \norm{d-z+y}{2}^2 ,
    \label{eq:L0formulation}
\end{equation}
where the L0 norm is replaced by its (convex) L1 counterpart. In this case the underlying model for $y$ is $y = \hat{y} + r + d$, where $r$ is zero-mean Gaussian noise and $d$ represents the ``oulier'' noise. Under this assumptions, $d$ practically becomes an ``indicator function'' of the presence of (sparse) outliers (see also~\cite{black91}). Besides the different angle in the derivation of the model, our approach differs from~\cite{mateos12} in two very important points. First, we use split-Bregman iteration by introducing the variable $b^k$ in the optimization procedure, see equations~(\refeq{eq:splitBregman1}) and~(\refeq{eq:splitBregman2}). In~\cite{mateos12} Equation~(\refeq{eq:L1formulationExtended}) is first solved using direct alternate minimization over $z$ and $d$, and then Equation~(\refeq{eq:L0formulation}) is solved via non-convex minimization using the previous solution as a starting point.
Second, considering the grid structure, we use the DCT approach to solve Equation~(\refeq{eq:splitBregman1_1}), instead of the classical Cholesky decomposition. The combination between split-Bregman and DCT results in a sound and fast algorithm for computing L1 splines on grid data.

\subsection{Handling missing data}

In the classical L2 formulation, a diagonal binary weighting matrix $W$ is used to cope with missing values (see Section~\ref{sec:robustL2} for details). Let us denote by $w$ the diagonal of $W$. Let us first define
\begin{align}
    \norm{z}{1, w} &= \sum_{\substack{i=1 \\ w_i = 1}}^{m} |z_i|
    & \text{and} &&
    \norm{z}{2, w} &= \left( \sum_{\substack{i=1 \\ w_i = 1}}^{m} z_i^2 \right)^{1/2} .
\end{align}
We then pose Equation~(\refeq{eq:WLSformulation}) as
\begin{gather}
    \hat{y} = \argmin{z} \norm{(z - y)}{2, w}^2 + s \norm{Dz}{2}^ 2 ,
\end{gather}
and equivalently extend Equation~(\refeq{eq:L1formulation}) as
\begin{gather}
    \hat{y} = \argmin{z} \norm{(z - y)}{1, w} + s \norm{Dz}{2}^ 2 .
\end{gather}
This leads to
\begin{equation}
    \min_{z, d} \ \norm{d}{1, w} + s \norm{Dz}{2}^ 2 + \tfrac{\lambda}{2} \norm{d - z + y}{2, w}^2 .
\end{equation}
Then the split-Bregman iteration can be written as
\begin{align}
    z^{k+1} &= \argmin{z} \, s \norm{Dz}{2}^ 2 + \tfrac{\lambda}{2} \norm{d^k - z + y - b^k}{2,w}^2 , \label{eq:WsplitBregman1_1} \\
    d^{k+1} &= \argmin{d} \, \norm{d}{1,w} + \tfrac{\lambda}{2} \norm{d - z^{k+1} + y - b^k}{2,w}^2 , \label{eq:WsplitBregman1_2} \\
    b^{k+1} &= b^k + (z^{k+1} - y - d^{k+1}) . \label{eq:WsplitBregman2}
\end{align}
Equation~(\refeq{eq:WsplitBregman1_1}) can be solved using Equation~(\refeq{eq:WLSsolutionDCT}). Solving Equation~(\refeq{eq:WsplitBregman1_2}) amounts to performing a shrinkage operation on the dimensions where $w$ equals 1. In Equation~(\refeq{eq:WsplitBregman2}), it suffices to update the dimensions of $b^k$ where $w$ equals 1.

\subsection{Handling multidimensional data}

Let us now return to the general case of $m$-dimensional data. Following Garcia~\cite{garcia10}, we extend Equation~(\refeq{eq:LSsolutionDCT}) as
\begin{equation}
    \hat{y} = \funct{DCT}^{-1}_m \left( \varGamma^m \circ \funct{DCT}_m (y) \right) ,
    \label{eq:LSsolutionDCTd}
\end{equation}
where $\funct{DCT}_m(\cdot)$ and $\funct{DCT}^{-1}_m(\cdot)$ stand for the $m$-dimensional DCT-II and inverse DCT-II functions, and $\circ$ denotes the Schur (element-wise) product. Notice that the multidimensional DCT is simply a composition of one-dimensional DCTs along each dimension. Extending Equation~(\refeq{eq:varGamma}), $\varGamma^m$ is  an $m$-th order tensor defined by
\begin{equation}
    \varGamma^m = 1^m  \div \left( 1^m + s \varLambda^m \circ \varLambda^m \right) ,
    \label{eq:varGammad}
\end{equation}
where $1^m$ is an $m$-th order tensor of ones, and $\div$ denotes the element-wise division.
Finally, $\varLambda^m$ is an $m$-th order tensor, defined by
\begin{equation}
    \varLambda^m_{i_1, \dots, i_m} = \sum_{j=1}^{d} \left(  -2 + 2 \cos \frac{(i_j - 1) \pi}{n_j} \right) .
    \label{eq:varLambdad}
\end{equation}
where $n_j$ denotes the size of $\varLambda^m$ along the $j$-th dimension.

\medskip

\noindent\textbf{The algorithm and its complexity.}
The pseudocode for the general $m$-dimensional case is presented in Algorithm~\ref{algo:l1spline}. Let us analyze its complexity.
The DCT and inverse DCT require $O(n \log n)$ operations, where $n=\prod_{1 \leq j \leq m} \, n_j$. The remaining operations are linear in $m$. The overall complexity of the algorithm is then $O \left( N_o N_i (m + n \log n) \right)$, where $N_o$ and $N_i$ are, respectively, the number of outer-loop and inner-loop iterations in Algorithm~\ref{algo:l1spline}. Notice that Goldstein and Osher~\cite{goldstein09} recommend to perform only one inner-loop iteration for achieving optimal efficiency. Thus, we set $N_i = 1$ for all experiments.
We will later see that in many cases the algorithm converges quickly ($N_o$ can be very small then).
The algorithm's complexity is thus dominated by the computation of the DCT and inverse DCT. Of course, these standard operations can be easily computed using GPU, speeding-up the execution by several orders of magnitude.

%
%

\begin{figure*}[t]
    \begin{boxedalgorithmic}
        \Function{L1Spline}{$y, s, \lambda, \varepsilon, N_i$}
        \State compute $\varGamma^m$ according to equations~(\refeq{eq:varGammad}) and~(\refeq{eq:varLambdad}).
        \State $d^1 \gets 0$, $b^1 \gets 0$, $k \gets 1$
        
        \Repeat
            \For{$i = 1$ \textbf{to} $N_i$} \Comment{usually $N_i = 1$.}
                \State $z^{k+1} \gets \funct{DCT}^{-1}_m (\varGamma^m \circ \funct{DCT}_m (d^k + y - b^k))$
                \Comment{$\circ$ denotes the element-wise product.}
                \State $\displaystyle d^{k+1} \gets \funct{Shrink} (z^{k+1} - y + b^k , 1/ \lambda)$
                \Comment{as defined in Equation~(\refeq{eq:shrinkage})}
            \EndFor
            \State $b^{k+1} \gets b^k + (z^{k+1} - y - d^{k+1})$
            \State $k \gets k+1$
        \Until{$\norm{z^{k+1} - z^k}{2} \, / \,  \norm{z^k}{2} > \varepsilon$}
        \State \Return $z^k$
        \EndFunction
    \end{boxedalgorithmic}

    \caption{Pseudocode of the L1 spline for evenly spaced-data via Bregman iteration}
    \label{algo:l1spline}
\end{figure*}

\section{Experimental results}
\label{sec:results}

For all experiments we adhere to the following setup:
\begin{enumerate}
    \item using generalized cross validation, we find the best estimate $\hat{s}$ for $s$ for the robust L2 formulation (problem~(\refeq{eq:WLSformulation}));
    \item we then find L2 splines (problem~(\refeq{eq:LSformulation})), robust L2 splines (problem~(\refeq{eq:WLSformulation})),\footnote{Code available at \url{http://www.mathworks.com/matlabcentral/fileexchange/25634-robust-spline-smoothing-for-1-d-to-n-d-data}.} and/or L1 splines (problem~(\refeq{eq:L1formulation})), setting $s = \hat{s}$.
\end{enumerate}
This protocol allows us to show that, even when $s$ is chosen to fit optimally the robust L2 formulation, the proposed method provides better estimates.
For the L1 formulation, in Equation~(\refeq{eq:L1formulation}), we simply set $\lambda=\min(s,1)$ for all examples. We recall that $N_i$ is set to 1. We also set $\varepsilon = 10^{-3}$ (see Algorithm~\ref{algo:l1spline}) and additionally limit the maximum number of outer iterations to a hundred. The algorithm stops when any of the two conditions is met.

\begin{figure*}
    \centerline{
        \hfill
        \subfloat[L2 spline]{
        \begin{minipage}{.32\textwidth}
                \includegraphics[width=\textwidth]{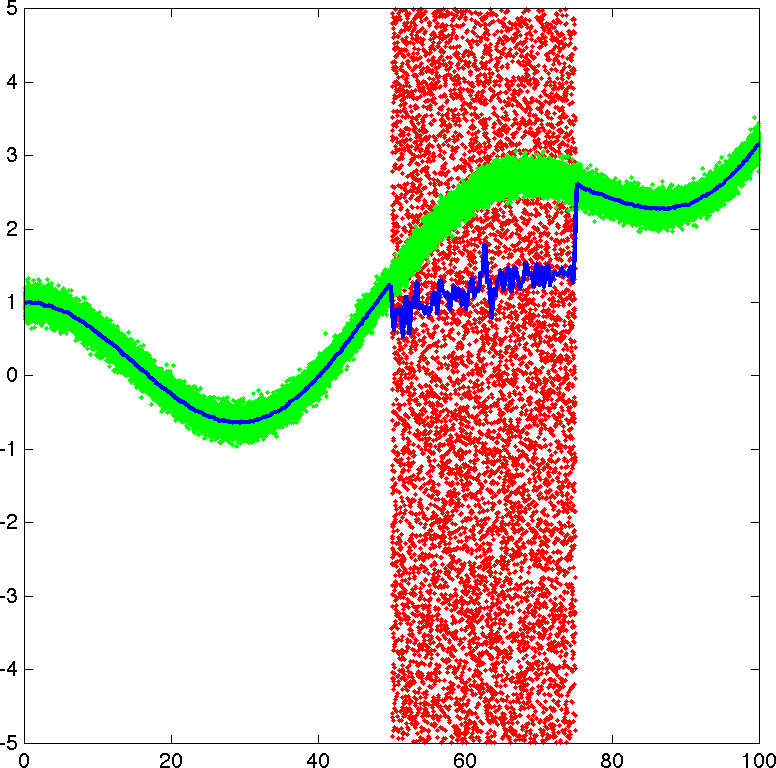} \\
                \includegraphics[width=\textwidth]{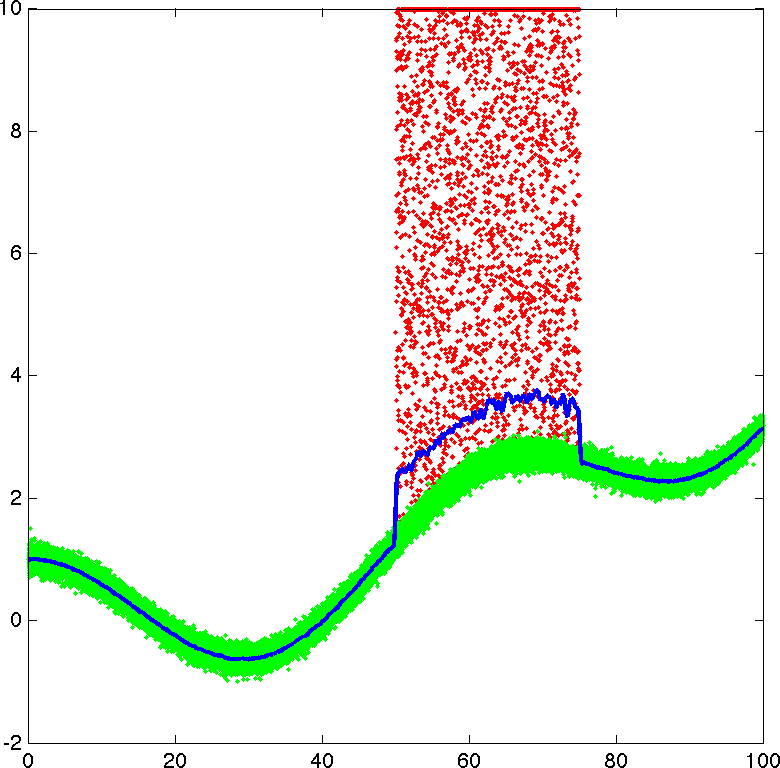}
        \end{minipage}
        }
        \hfill
        \subfloat[Robust L2 spline]{
        \begin{minipage}{.32\textwidth}
                \includegraphics[width=\textwidth]{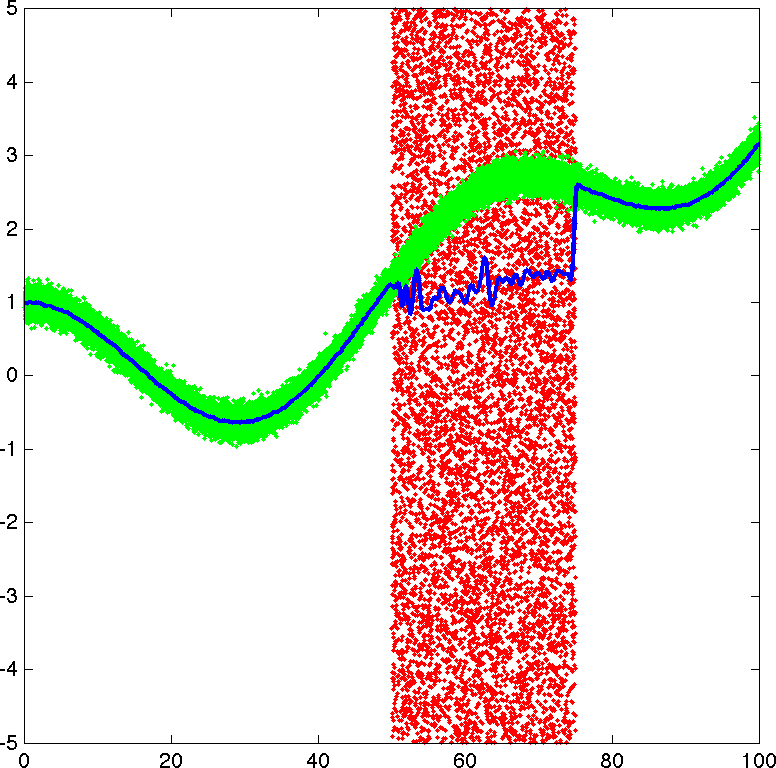} \\
                \includegraphics[width=\textwidth]{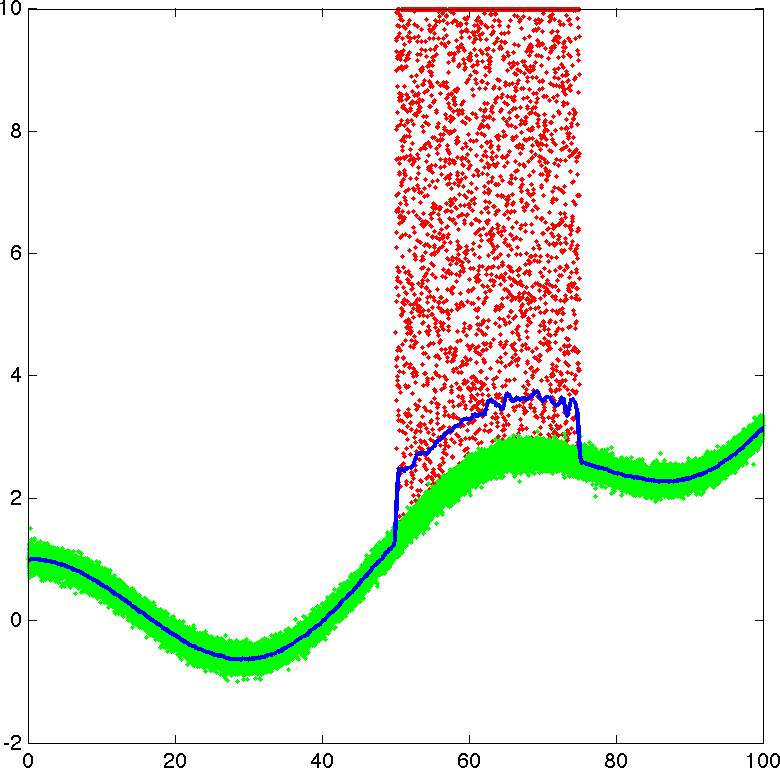}
        \end{minipage}
        }
        \hfill
        \subfloat[L1 spline]{
        \begin{minipage}{.32\textwidth}
                \includegraphics[width=\textwidth]{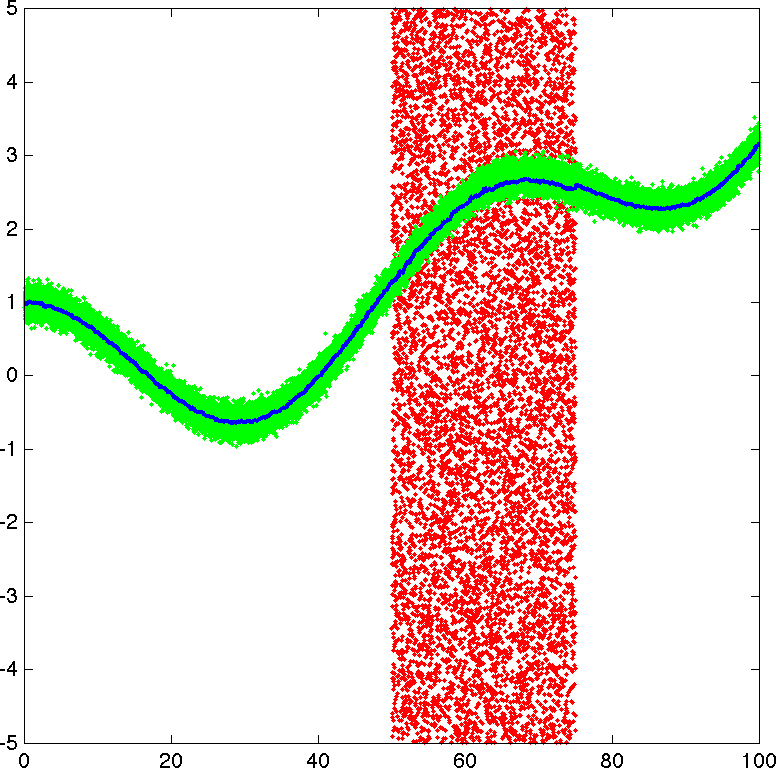} \\
                \includegraphics[width=\textwidth]{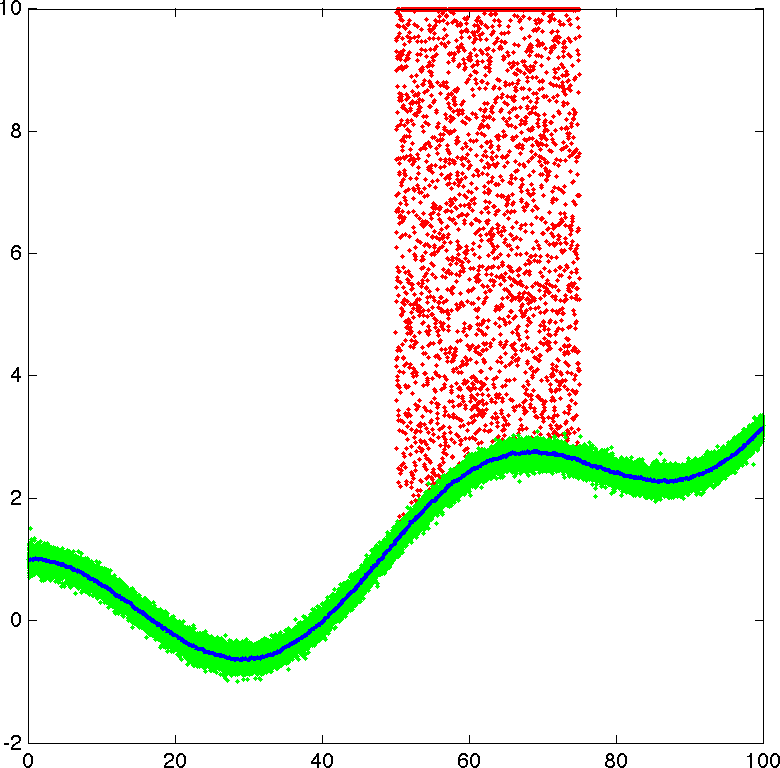}
        \end{minipage}
        }
        \hfill
    }
    \caption{Noisy data (in green) further contaminated with uniform noise (in red). The spline $\hat{y}$, depicted in blue, is obtained using different fitting terms. The proposed method with the L1 fitting term recovers the correct shape.}
    \label{fig:spline_1D_16}
\end{figure*}

Fig.~\ref{fig:spline_1D_16} presents two one-dimensional examples. We depict the original signal $\hat{y} \in [1, \dots, n] \rightarrow \Real$, where $n=2^{16}$. We observe the signal $y = \hat{y} + r_1$, where $r_1$ is Gaussian noise.
Some points $y_j$ ($1 \leq j \leq n$) are further contaminated with uniform noise $r_2$, where $r_2 \in [ a, \dots, b ]$, such that $y_j = \min (\max (\hat{y}_j + r_1 + r_2, a), b)$. The points affected by $r_2$ are depicted in red and the remaining ones in green.
In the top row, $a = -5, b = 5$; and in the bottom row, $a = 0, b = 5$. In both examples, only the L1 spline is correct. The classical and robust L2 splines are both unable to correctly recover the original data in the corrupted part.

\begin{figure}
    \centering
    \includegraphics[width=.5\columnwidth]{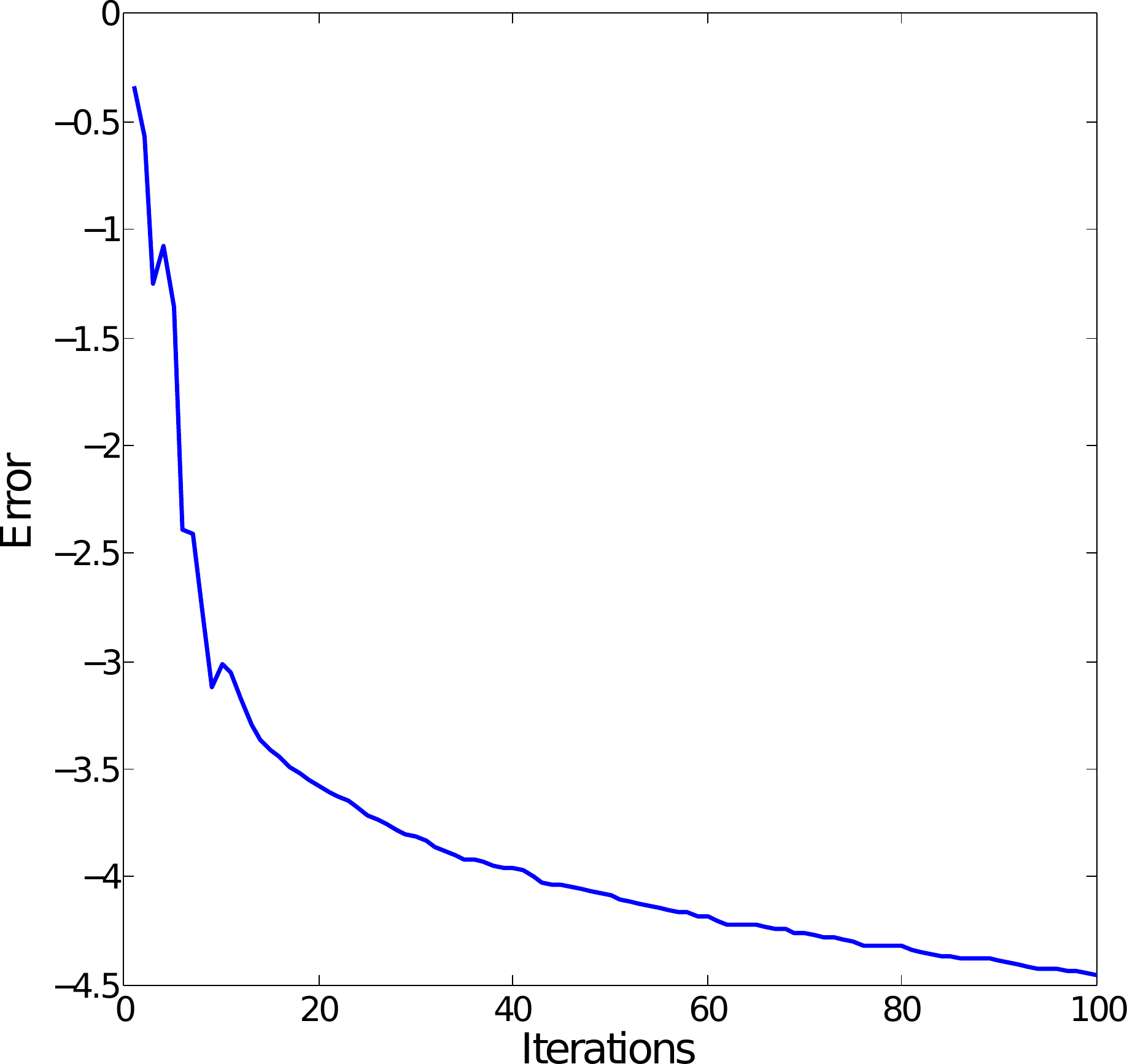}

    \caption{The relative error (in logarithmic scale) as the iterations progress when computing the first example in Figure~\protect\ref{fig:spline_1D_16}. The algorithm is able to quickly decrease the error during the first 20 iterations.}
    \label{fig:errorPlot}
\end{figure}

Fig.~\ref{fig:errorPlot} shows the evolution of the relative error as the number of iterations increases for the example in Fig.~\ref{fig:spline_1D_16} (top row).
As we can observe, the proposed algorithm is able to converge quickly, reaching a precision of $10^{-3}$ in less than 20 iterations.

\begin{figure}
    \centering
    \includegraphics[width=.5\columnwidth]{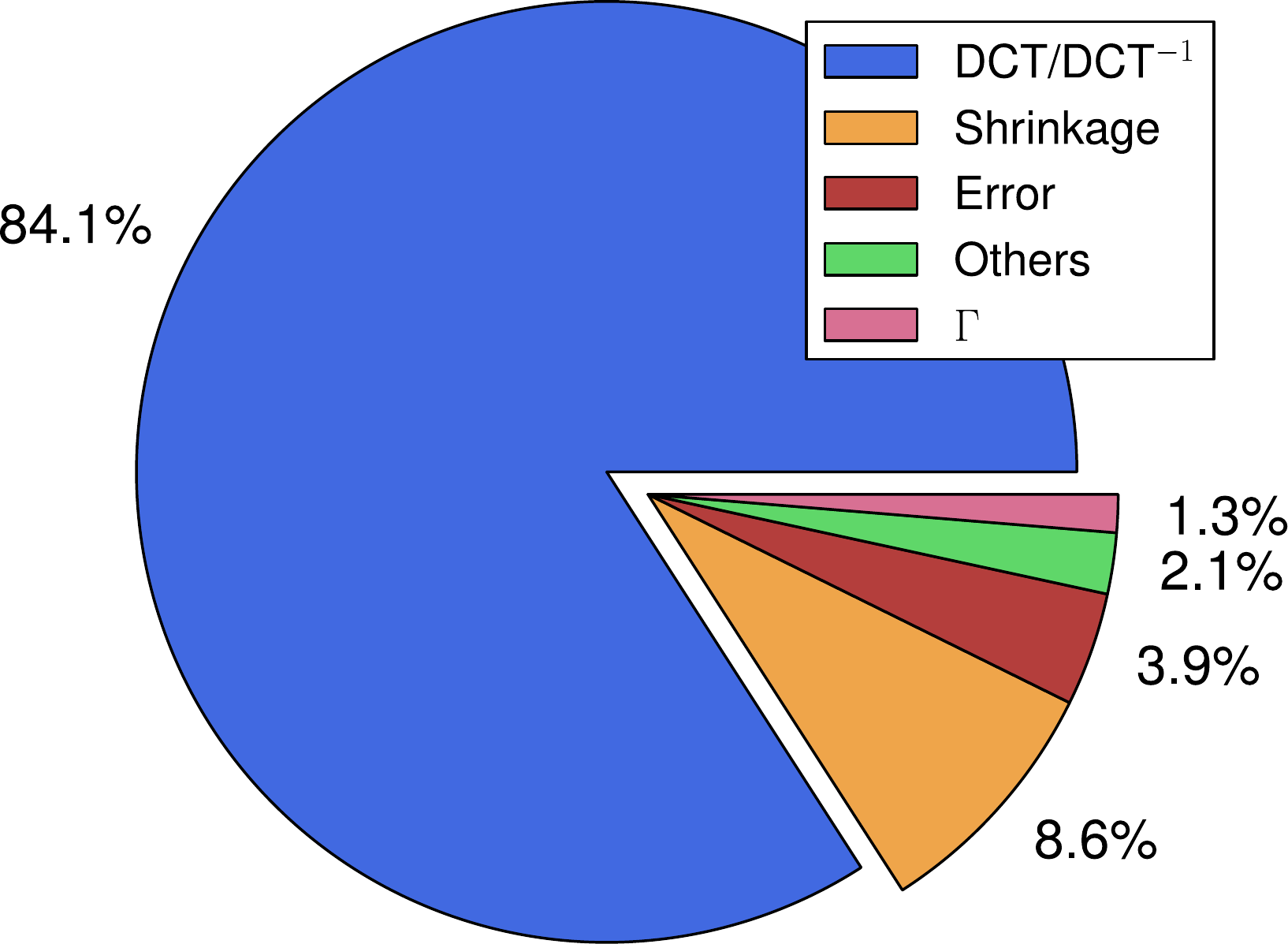}

    \caption{Percentage of the execution time spent in each operation when computing the first example in Figure~\protect\ref{fig:spline_1D_16}. Clearly, the vast majority of time is spent in DCT or inverse DCT operations.}
    \label{fig:timeDistribution}
\end{figure}

Fig.~\ref{fig:timeDistribution} depicts the relative time-cost of each operation during the execution of the proposed method (Fig.~\ref{fig:spline_1D_16}, top row). Computing the DCT and the inverse DCT covers more than 84\% of the total running time.
Implementing these standard operations in GPU would boost the performance of the algorithm by orders of magnitude.

\begin{figure*}
    \centering
    \begin{tabular}{@{\hspace{4pt}}c@{\hspace{4pt}}c@{\hspace{4pt}}c@{\hspace{4pt}}}
        \subfloat[Original data $y$]{
            \includegraphics[width=.3\textwidth]{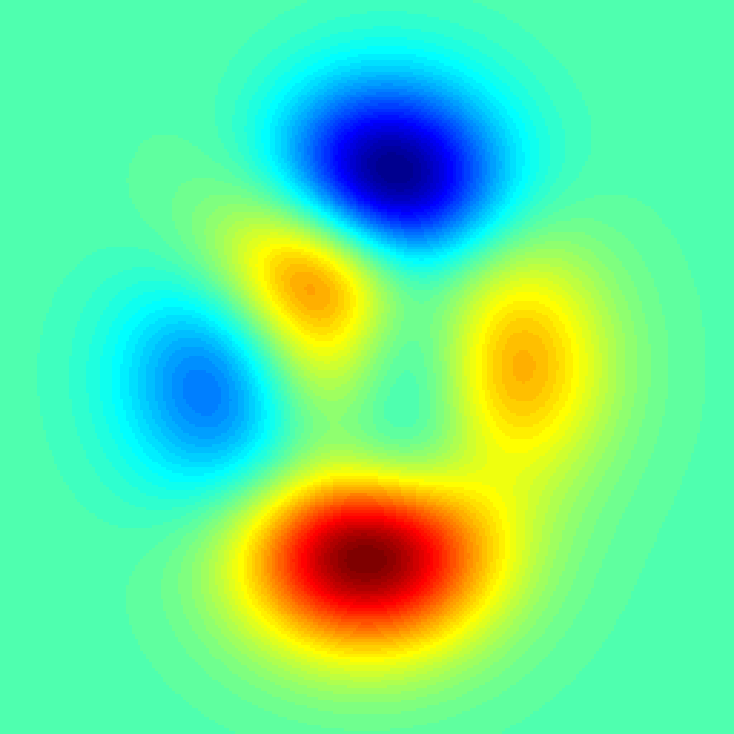}
            \label{fig:spline_2D_n30_original}
        }
        &
        \subfloat[Noisy data $y+r_1$]{
            \includegraphics[width=.3\textwidth]{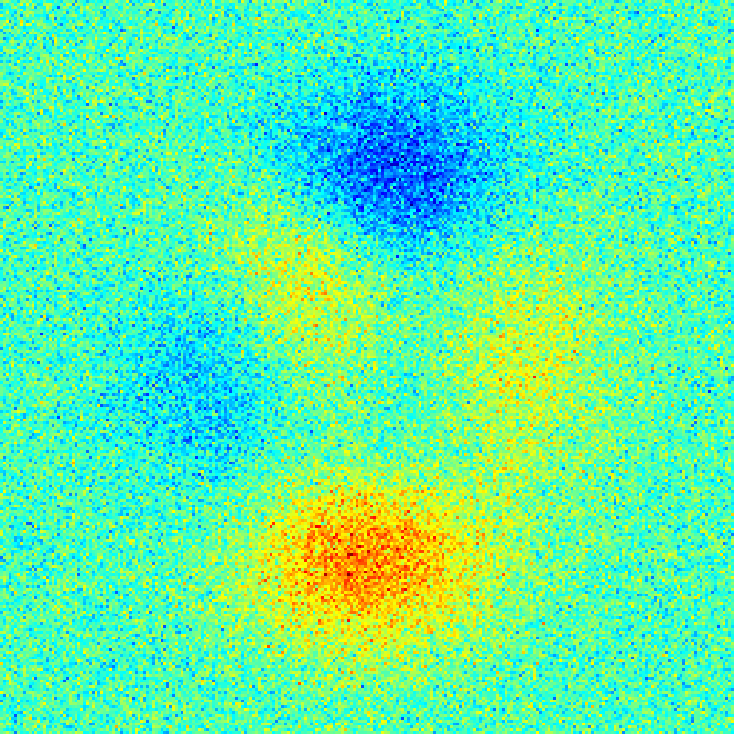}
            \label{fig:spline_2D_n30_noisy}
        }
        &
        \subfloat[Corrupted data $y+r_1+r_2$]{
            \includegraphics[width=.3\textwidth]{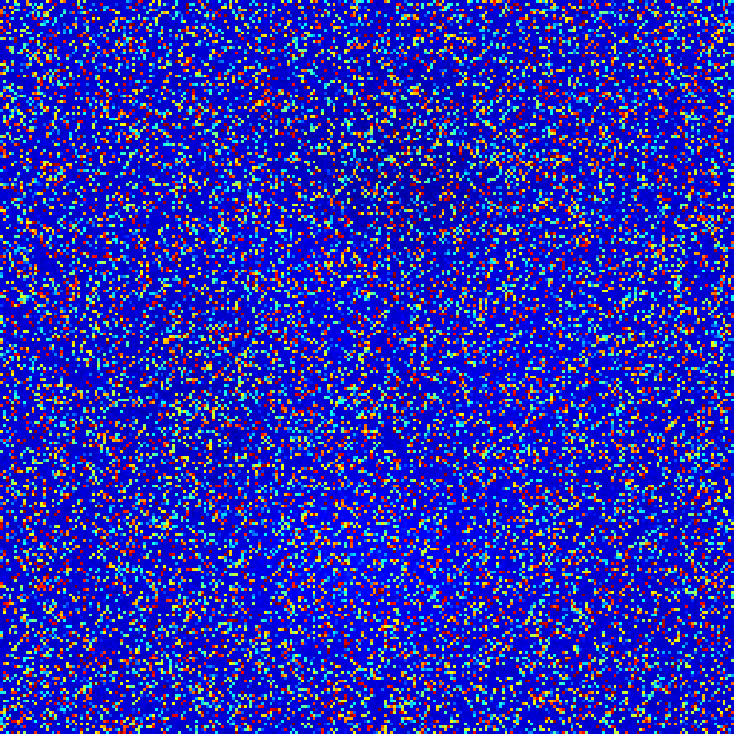}
            \label{fig:spline_2D_n30_corrupted}
        }
        \tabularnewline

        \subfloat[L2 spline]{
            \includegraphics[width=.3\textwidth]{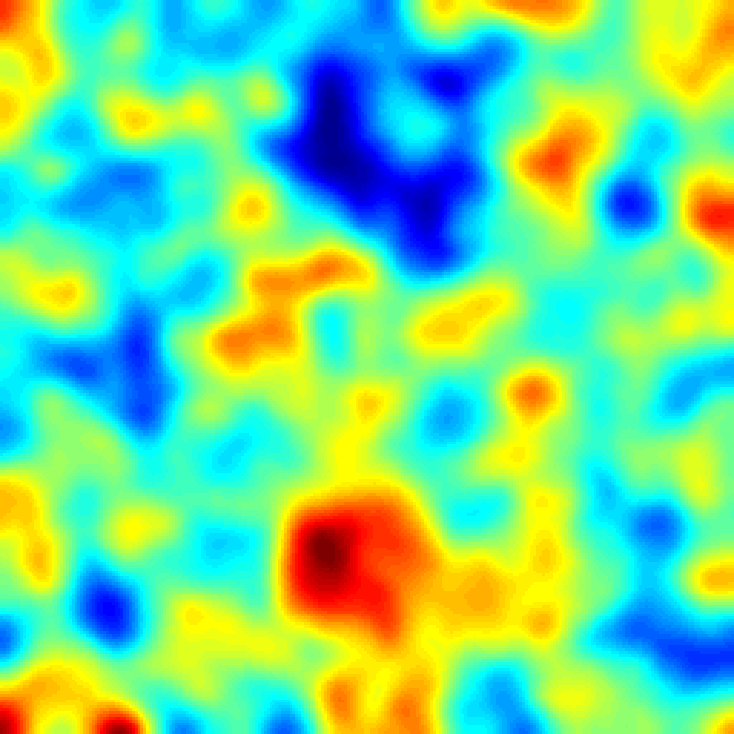}
        }
        &
        \subfloat[Robust L2 spline]{
            \includegraphics[width=.3\textwidth]{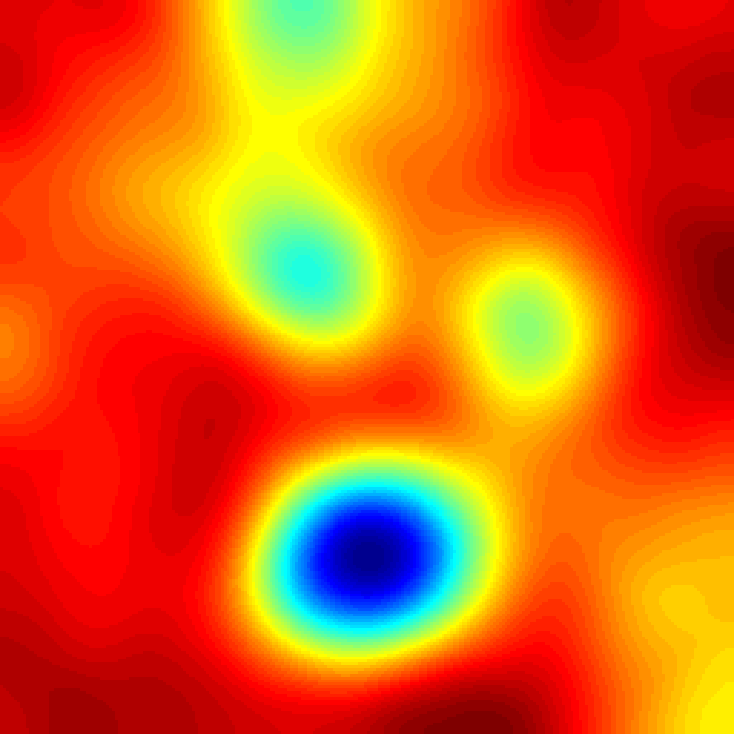}
        }
        &
        \subfloat[L1 spline]{
            \includegraphics[width=.3\textwidth]{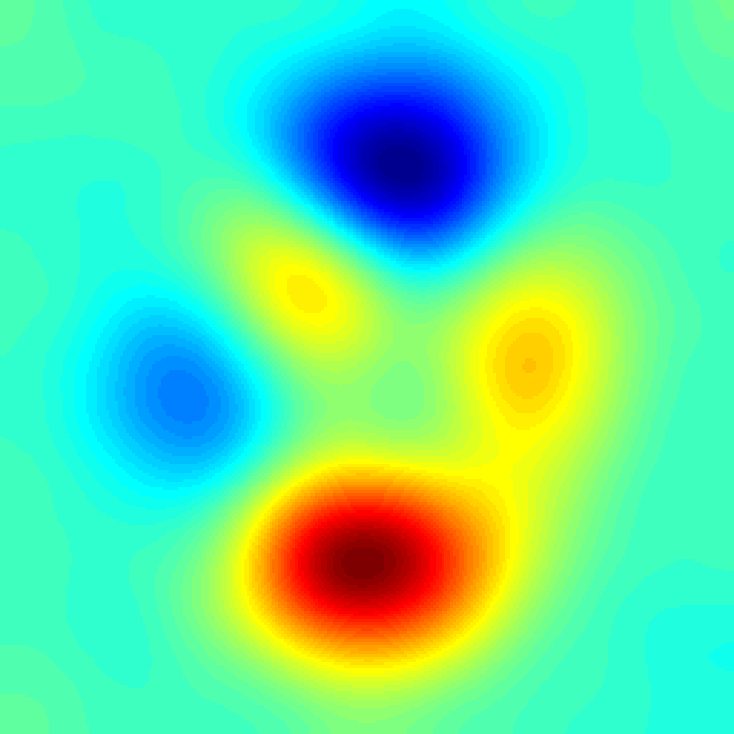}
        }
        \tabularnewline
    \end{tabular}
    \caption{Synthetic two-dimensional example: the original data $y$ is contaminated with Gaussian noise $r_1$ and then with a large uniform noise $r_2$. With this input signal, $y+r_1+r_2$, the proposed method is the only one able to recover the correct shape.}
    \label{fig:spline_2D_n30}
\end{figure*}

In Fig.~\ref{fig:spline_2D_n30} we present a two-dimensional example. We depict the original signal $\hat{y} \in [1, \dots, 256]^2\rightarrow [-6.5497, 8.1054]$ in Fig.~\ref{fig:spline_2D_n30_original}, and we add two types of noise: first, Gaussian noise $r_1$ with zero mean and variance $\sigma^2=2$ (Fig.~\ref{fig:spline_2D_n30_noisy}), and then uniform noise $r_2$ in the interval $[ -5 \cdot \max (\hat{y} + r_1), \dots, 5 \cdot \max (\hat{y} + r_1) ]$ (Fig.~\ref{fig:spline_2D_n30_corrupted}). Again, only the L1 spline correctly recovers the original signal.

\begin{figure}
    \centering{
        \includegraphics[width=.7\columnwidth]{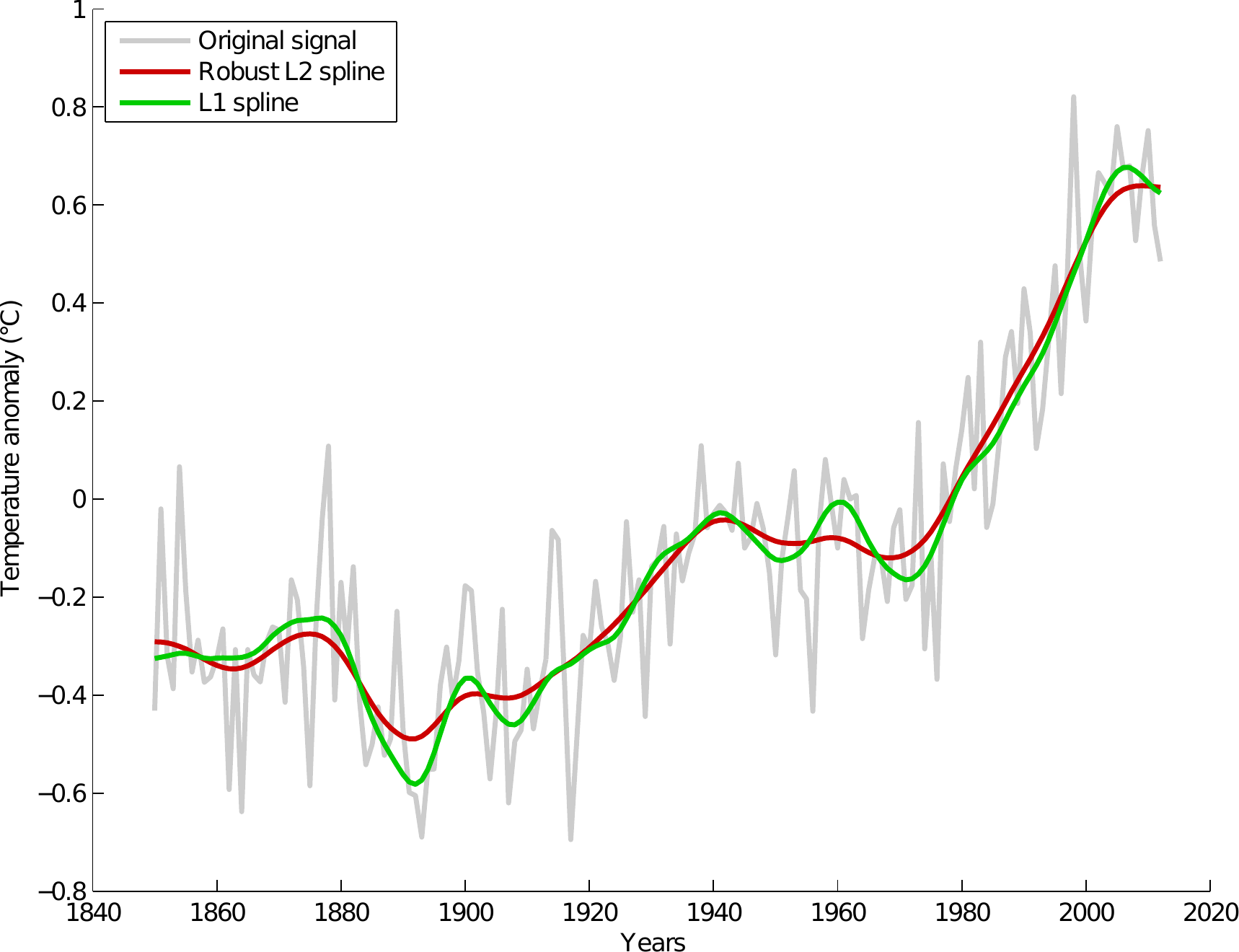}
    }
    
    \caption{Global average land temperature anomaly (\textcelsius) with respect to 1961-1990: smoothed versus original year-averaged data~\cite{brohan06}. Control points are points such that $\hat{y}_i = y_i$.}
    \label{fig:temp}
\end{figure}

We next test the proposed algorithm with a climate time-series provided by the Met Office Hadley Centre~\cite{brohan06}.\footnote{Data are available in \url{http://hadobs.metoffice.com/crutem3/diagnostics/global/nh+sh/annual}.} The dataset contains the evolution of global average land temperature anomaly (in \textcelsius) with respect to the 1961-1990 average temperature. The results, which confirm an upward trend in the second half of the 20th century, are shown in Fig.~\ref{fig:temp}.

\begin{figure*}
    \centering
        \subfloat[Variable $s$, fixed $\lambda=1$.]{
            \includegraphics[width=.7\columnwidth]{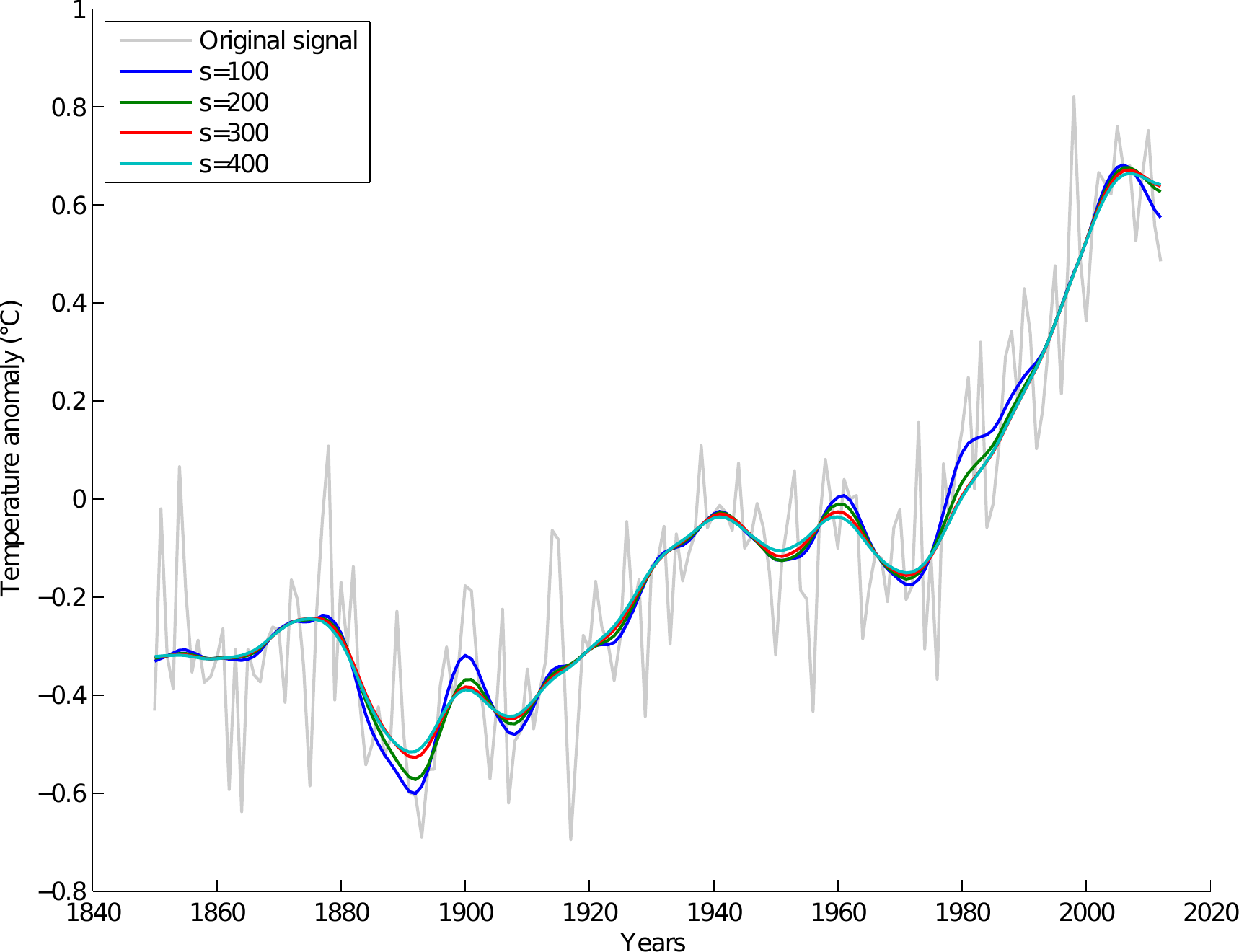}
            \label{fig:params_s}
        }

        \subfloat[Fixed $s=200$, variable $\lambda$.]{
            \includegraphics[width=.7\columnwidth]{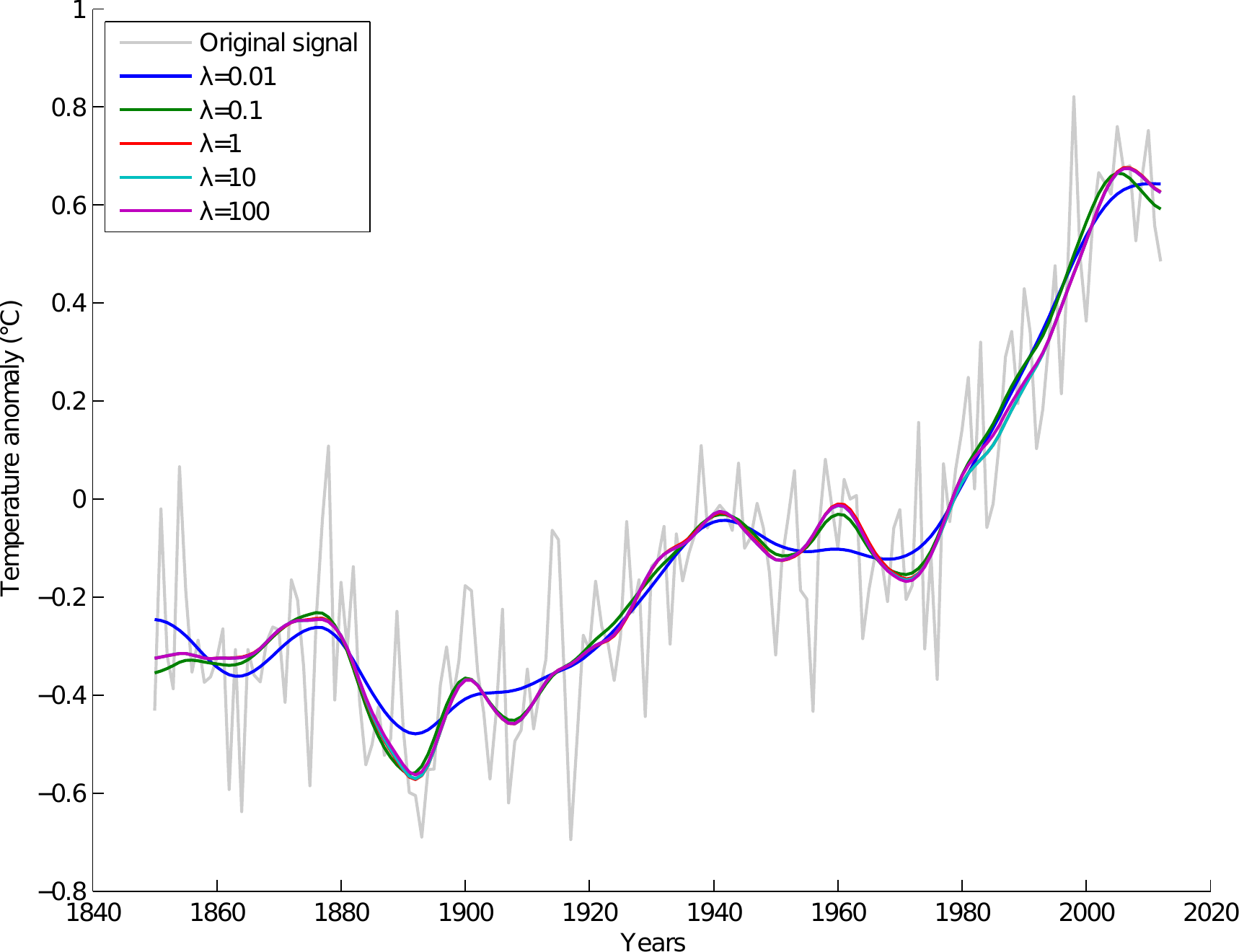}
            \label{fig:params_lambda}
        }
    
    \caption{Changing the value of $s$ in L1 splines controls the trade-off between fitting and smoothing. L1 splines are quite insensitive to the value of $\lambda$, making it easy to tune.}
    \label{fig:params}
\end{figure*}

We also test on this dataset the effect of varying the parameters $s$ and $\lambda$, see Fig.~\ref{fig:params}. As in the classical L2 formulation, $s$ has a direct impact on the obtained result, controlling the degree of smoothness of the solution, see Fig.~\ref{fig:params_s}. On the contrary, Fig.~\ref{fig:params_lambda} shows that the newly introduced parameter $\lambda$ is very stable and provides very similar results in a wide range $(\lambda \in [0.1, 100])$. This stability allows us to fix its value to $\lambda=1$ for all the experiments in this work.

\begin{figure}
    \centering{
        \includegraphics[width=.7\columnwidth]{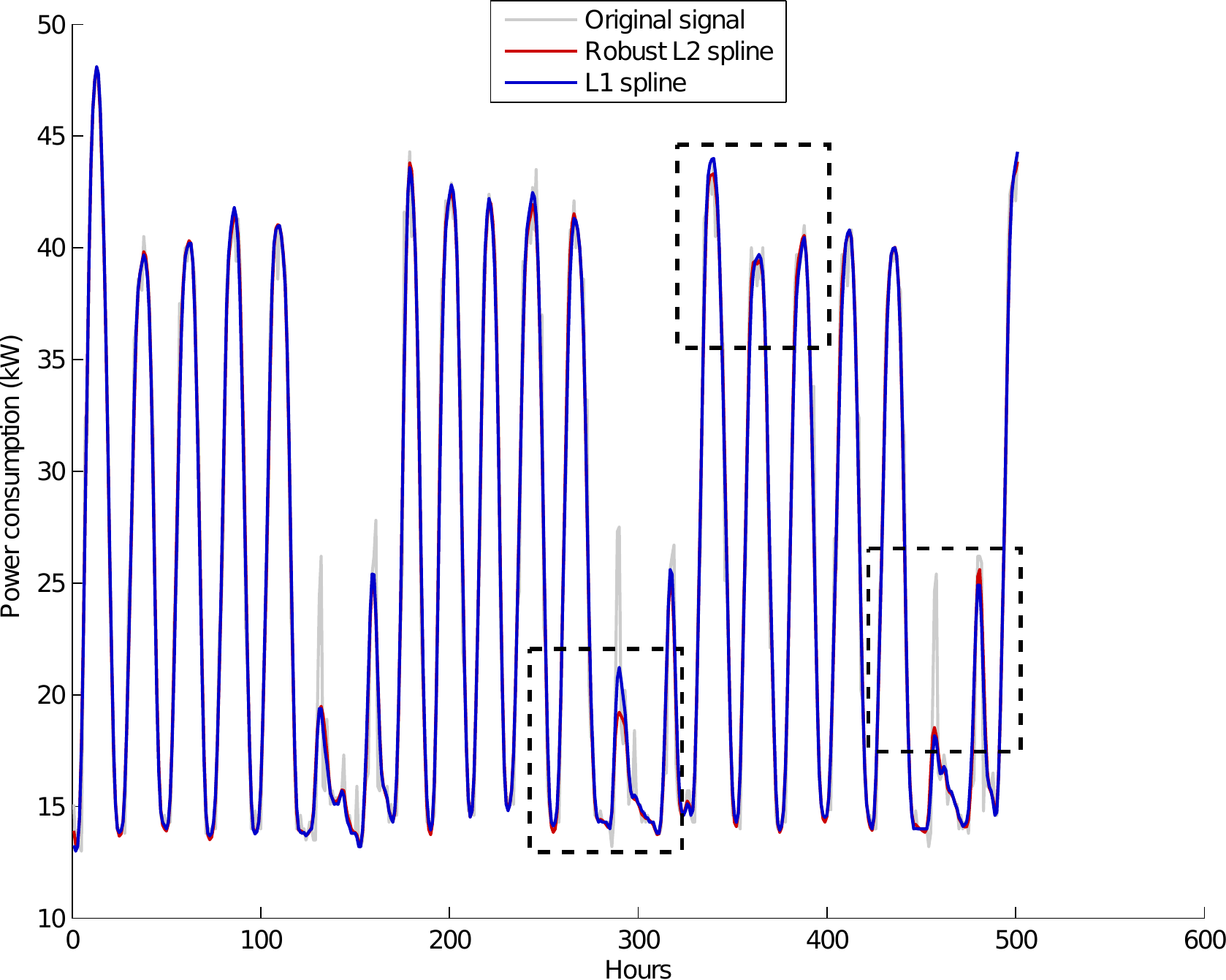}
    }
    
    \centerline{
        \hfill
        \includegraphics[width=.2\columnwidth]{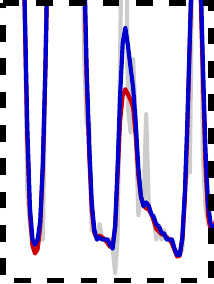}
        \hfill
        \includegraphics[width=.2\columnwidth]{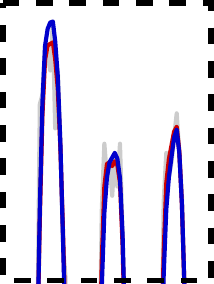}
        \hfill
        \includegraphics[width=.2\columnwidth]{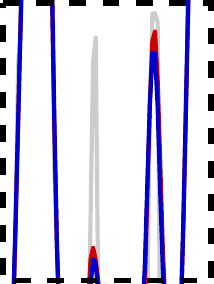}
        \hfill
    }
    \caption{Power consumption measurements (in kW) for a government building~\protect\cite{mateos12} with zoom-in details on the bottom.}
    \label{fig:load}
\end{figure}

Another interesting example is presented in~\cite{mateos12}. The dataset consists of power consumption measurements (in kW) for a government building, collected every fifteen minutes from July 2005 to October 2010. As in~\cite{mateos12}, we downsample the data by a factor of four, yielding one measurement per hour, and use only a subset of the whole data. The results are displayed in Fig.~\ref{fig:load}.

\begin{figure*}
    \centerline{
    \hfill
    \subfloat[L2 spline]{
        \includegraphics[width=0.3\textwidth]{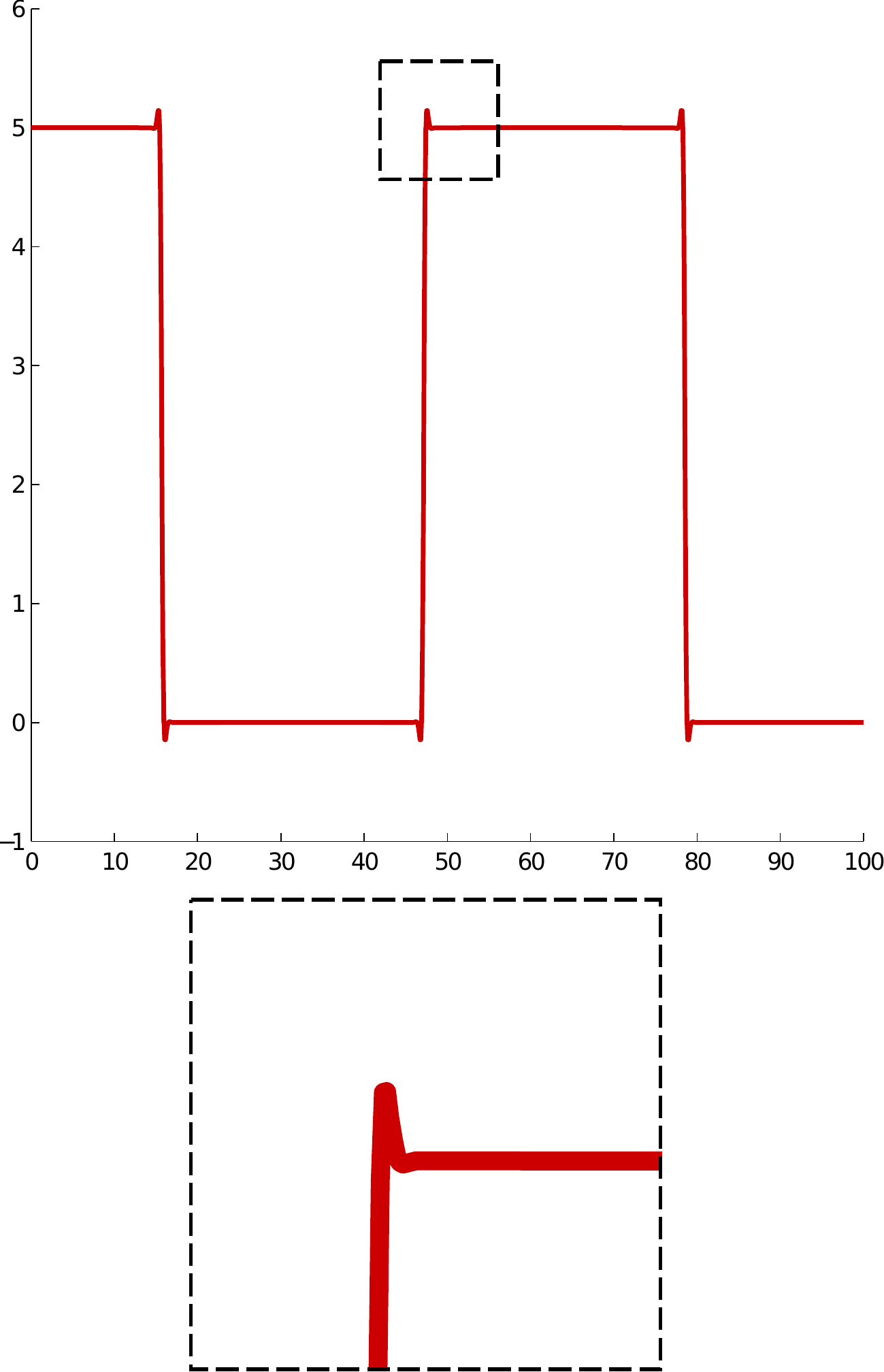}
    }
    \hfill
    \subfloat[Robust L2 spline]{
        \includegraphics[width=0.3\textwidth]{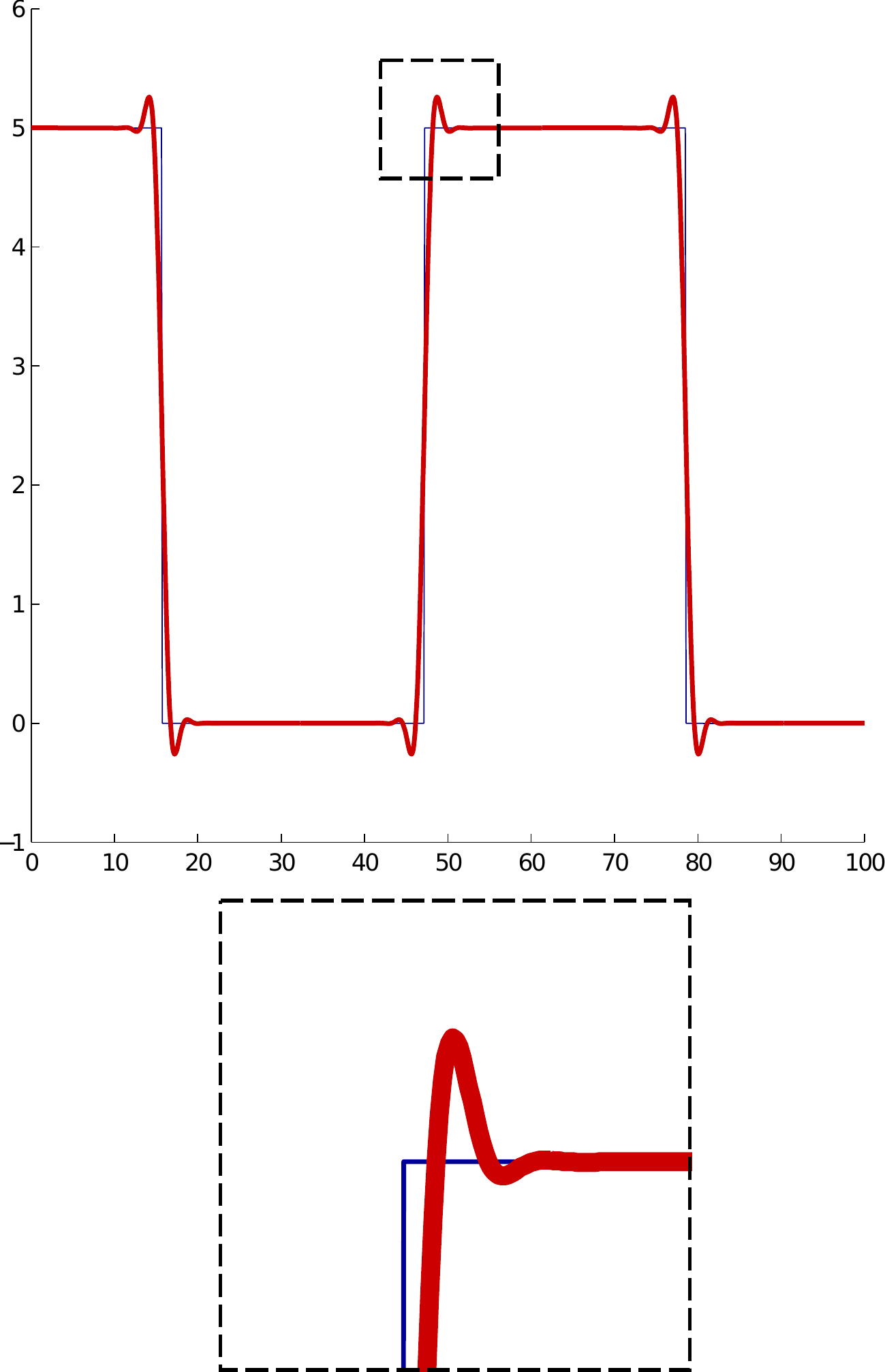}
    }
    \hfill
    \subfloat[L1 spline]{
        \includegraphics[width=0.3\textwidth]{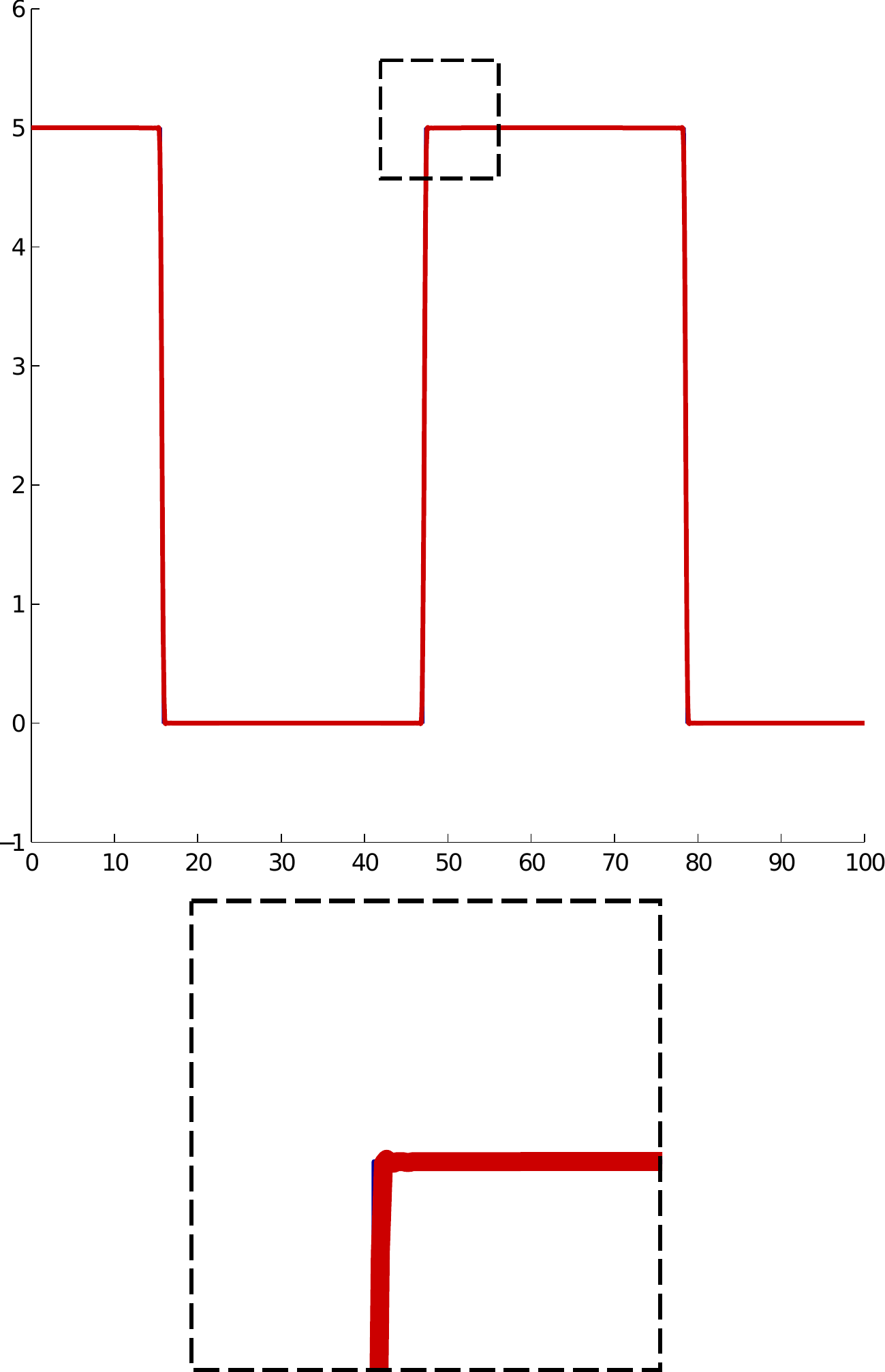}
    }
    \hfill
    }
    
    \caption{In the case of a periodic piece-wise constant signal, the L2 and robust L2 approximations exhibit overshoot and ringing (see zoom-in details on the bottom). These effects are much attenuated by the L1 spline approximation.}
    \label{fig:square}
\end{figure*}

We also test in a synthetic example the ability to recover signals with sharp transitions, see Fig.~\ref{fig:square}. In this case we use a simple piece-wise constant function. We can observe clear overshoot (plus ringing) effects on the L2 and robust L2 splines. The robust L2 spline also results in transitions with less vertical slopes, creating a bluring effect. With the L1 spline we obtain a much better reconstruction, with almost non-existent overshooting.

This very same effect can be observed in real examples, see Fig.~\ref{fig:119082}. When approximating images with splines, some structure is lost by blur and some structure is artificially created by overshooting and ringing. This can be observed in Figs.~\ref{fig:119082_original} and~\ref{fig:119082_noise}, were the difference between the original image and the image estimated by robust L2 splines exhibits structure. Observe, however, that almost no structure in the difference is visible when the reconstruction is performed using L1 splines.

\newlength{\mylength}
\setlength{\mylength}{0.234\textwidth}
\setlength{\fboxsep}{0pt}
\setlength{\fboxrule}{1pt}

\begin{figure*}
    \begin{tabular}{@{\hspace{0pt}}c@{\hspace{0pt}}c@{\hspace{0pt}}}
        \subfloat[]{
            \begin{minipage}{.5\columnwidth}
            \begin{tabular}{@{\hspace{0pt}}m{10pt}@{\hspace{0pt}}m{\mylength}@{\hspace{2pt}}m{\mylength}@{\hspace{0pt}}}
            
                & \multicolumn{2}{c}{\begin{footnotesize}Original image\end{footnotesize}} \\
                & \multicolumn{2}{c}{\includegraphics[width=\mylength]{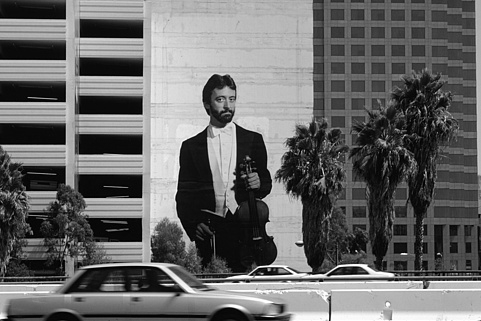}} \\
                
                &
                \multicolumn{1}{c}{\begin{footnotesize}Robust L2 spline\end{footnotesize}} &
                \multicolumn{1}{c}{\begin{footnotesize}L1 spline\end{footnotesize}} \\
                
                \begin{sideways}\begin{footnotesize}Approximation\end{footnotesize}\end{sideways} &
                \includegraphics[width=\mylength]{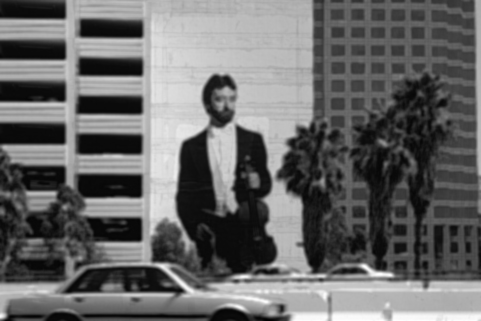} &
                \includegraphics[width=\mylength]{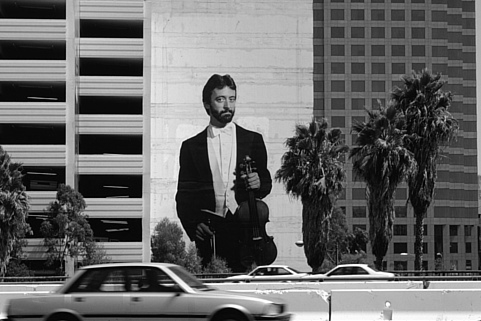} \\
                
                \begin{sideways}\begin{footnotesize}Difference\end{footnotesize}\end{sideways} &
                \fbox{\includegraphics[width=\mylength]{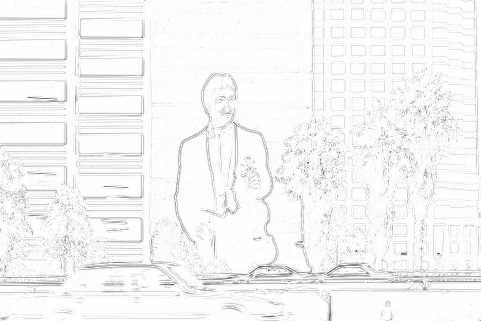}} &
                \fbox{\includegraphics[width=\mylength]{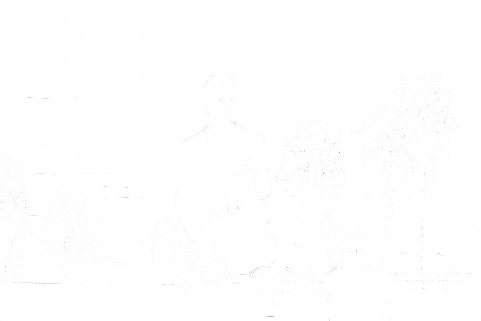}} \\
            \end{tabular}
            \end{minipage}
            \label{fig:119082_original}
        }
        \subfloat[]{
            \begin{minipage}{.5\columnwidth}
            \begin{tabular}{@{\hspace{0pt}}m{\mylength}@{\hspace{2pt}}m{\mylength}@{\hspace{0pt}}}

                \multicolumn{2}{c}{\begin{footnotesize}Noisy image\end{footnotesize}} \\
                \multicolumn{2}{c}{\includegraphics[width=\mylength]{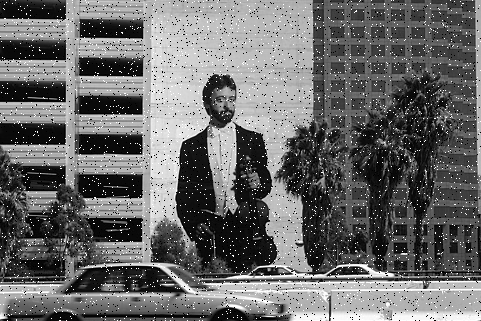}} \\
                
                \multicolumn{1}{c}{\begin{footnotesize}Robust L2 spline\end{footnotesize}} &
                \multicolumn{1}{c}{\begin{footnotesize}L1 spline\end{footnotesize}} \\
                
                \includegraphics[width=\mylength]{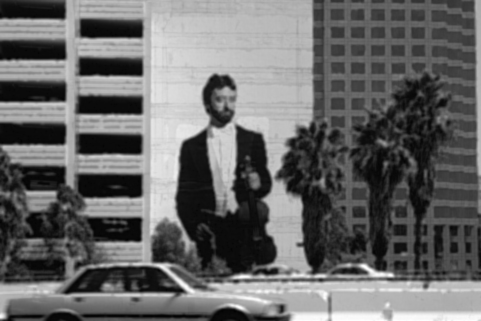} &
                \includegraphics[width=\mylength]{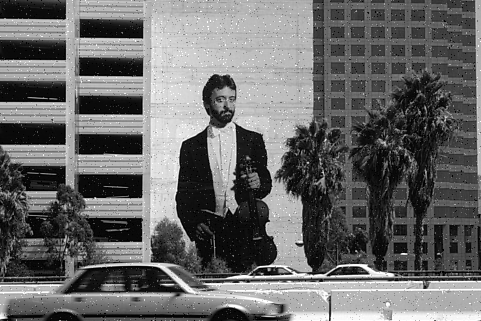} \\
                
                \fbox{\includegraphics[width=\mylength]{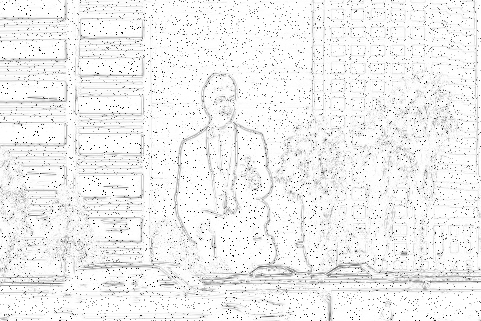}} &
                \fbox{\includegraphics[width=\mylength]{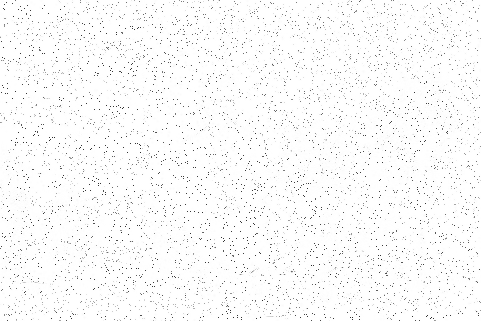}} \\
            \end{tabular}
            \end{minipage}
            \label{fig:119082_noise}
        }
    \end{tabular}
    
    \caption{\protect\subref{fig:119082_original} When no noise is added, the robust L2 spline removes much more structure than the L1 spline. For the robust L2 spline and for the L1 spline, the norm of the difference is 9.9866 and 0.7301, respectively. \protect\subref{fig:119082_noise} When salt \& pepper noise is added, the robust L2 spline removes noise and structure, while the L1 spline removes much less structure but leave some noise. Notice that we are not proposing a new image-denoising algorithm, we just use images for showing that L1 splines respect more the signal structure than L2 splines.}
    \label{fig:119082}

\end{figure*}

\subsection{Application to range data}

In this section we perform smoothing of depth data obtained with a Kinect camera. This kind of data is particularly challenging because:
\begin{itemize}
    \item it presents relatively smooth areas separated by sharp transitions,
    \item edges are highly noisy, that is, edge pixels oscillate over time between foreground and background, and
    \item it contains missing data, which appear for two different reasons: (1) the disparity between the IR projector and the IR camera produces ``shadows,'' and (2) the depth cannot be recovered in areas where the IR pattern is not clearly observable (e.g., because they receive direct sunlight or interference from another Kinect).
\end{itemize}
We use splines to interpolate and denoise these data, showing the advantage of L1 splines over its robust L2 counterpart.
The displayed images are part of the LIRIS human activities dataset~\cite{harl2012}.

\begin{figure*}
    \centering
    \begin{tabular}{ccc}
        \subfloat[Original depth image]{
            \includegraphics[width=\mylength]{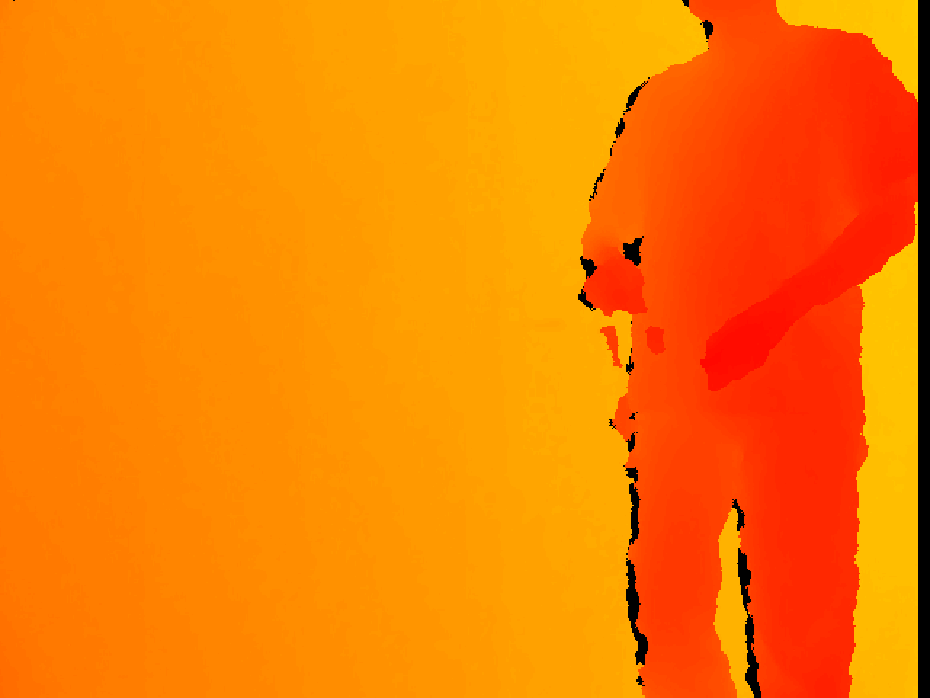}
            \label{fig:kinect2D_original}
        }
        \subfloat[Robust L2 spline]{
            \includegraphics[width=\mylength]{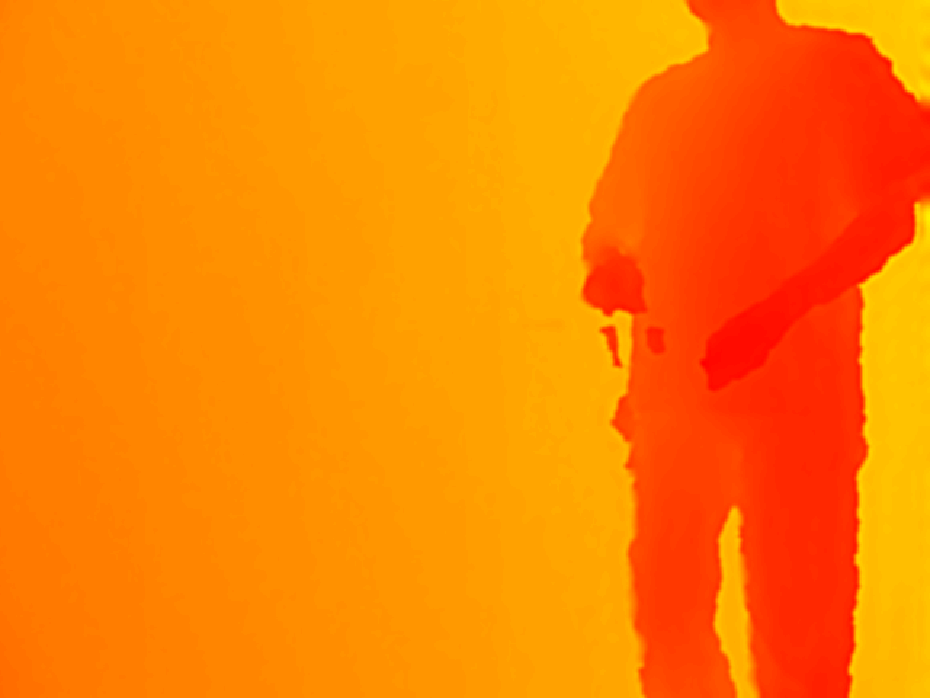}
        }
        \subfloat[L1 spline]{
            \includegraphics[width=\mylength]{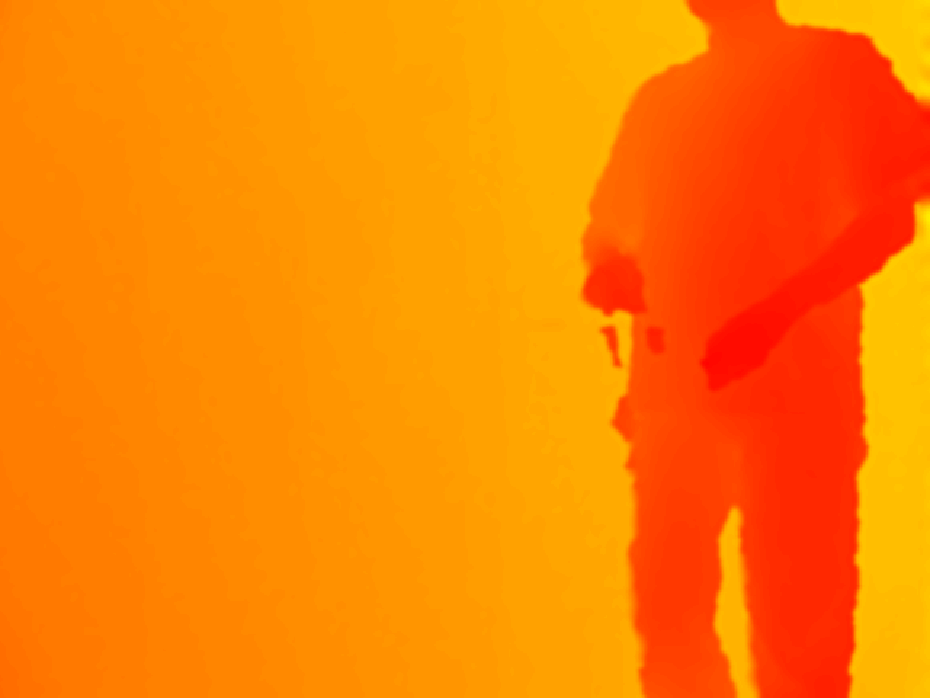}
        }
    \end{tabular}

    \subfloat[Profile of a single row. Left, original values; right, original and reconstructed values.]{
        \centerline{
            \hfill
            \includegraphics[width=0.45\textwidth]{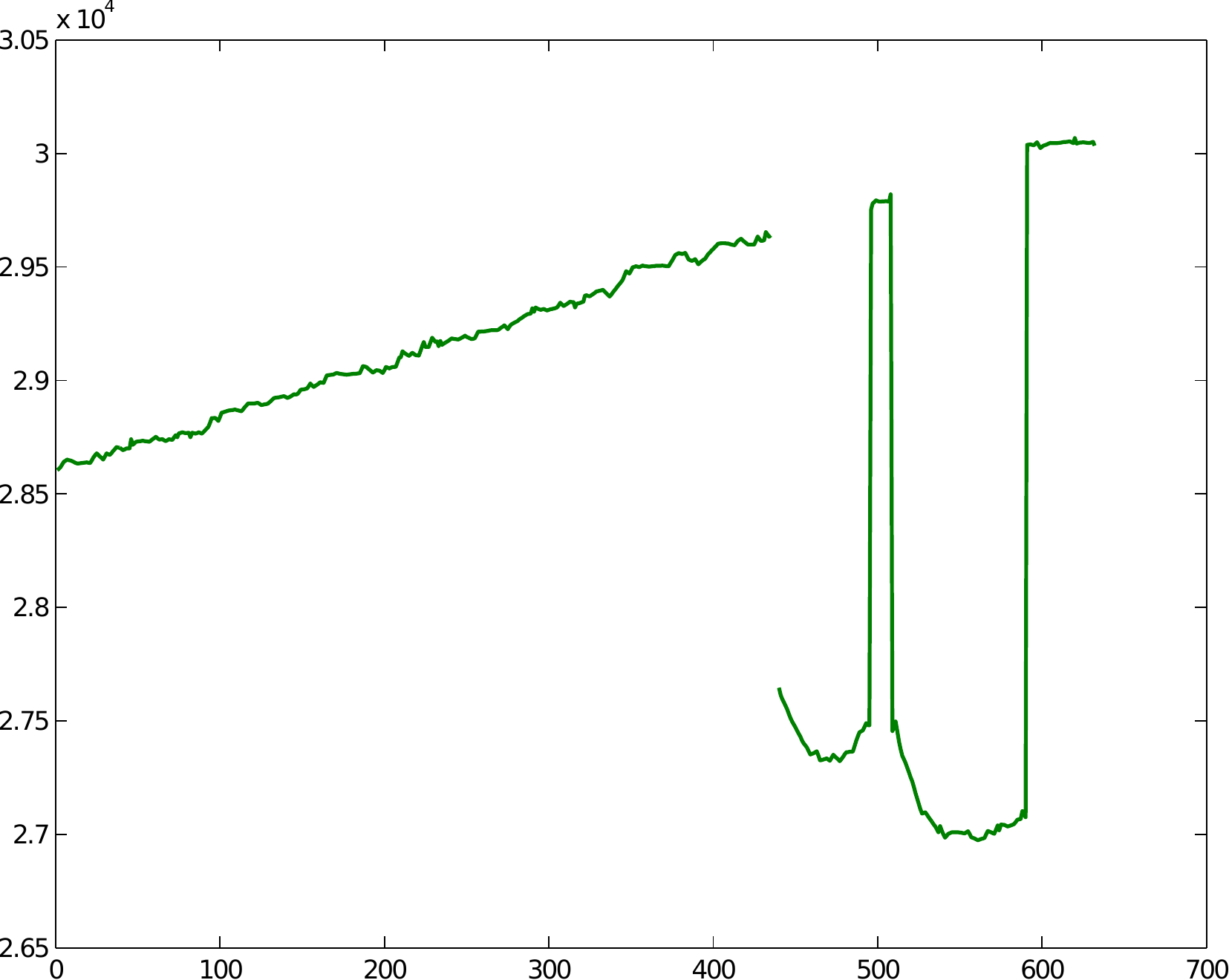}
            \hfill
            \includegraphics[width=0.45\textwidth]{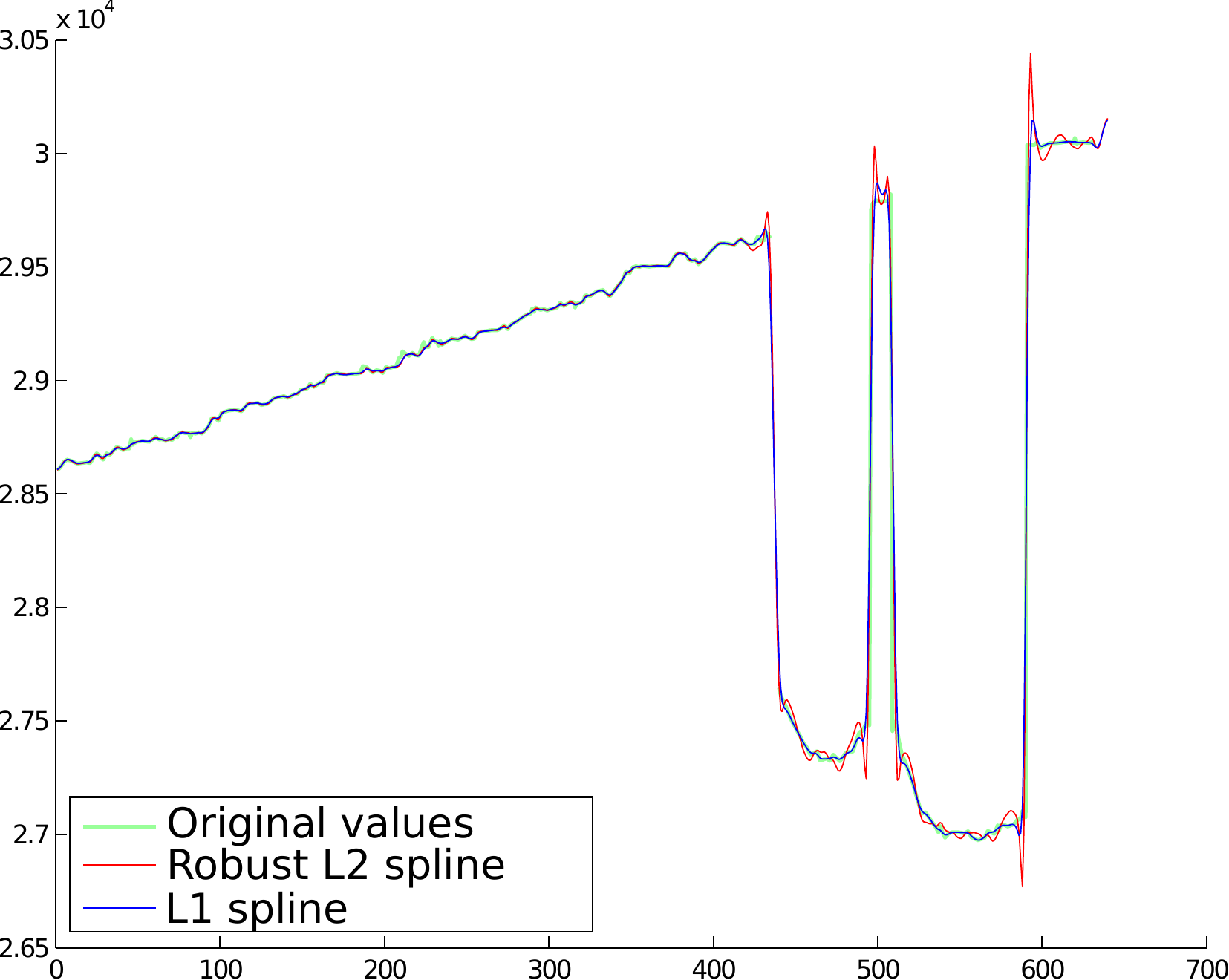}
            \hfill
        }
    }
    \caption{Smoothing and interpolating using a single depth frame. We simply fit a 2D spline to the depth image. Since the data presents several ``jumps,'' the robust L2 spline must under-smooth the data to be able to fit it correctly. The L1 spline presents a good trade-off between fitting and smoothing.}
    \label{fig:kinect2D}
\end{figure*}

In the first example, shown in Fig.~\ref{fig:kinect2D}, we use a single depth frame (with standard Kinect resolution of $640 \times 480$). The missing data are represented in black, while depth data goes from red to yellow as depth increases. Both, the L1 spline and the robust L2 spline are able to interpolate the missing data with reasonable values. Notice, however, that the latter exhibits, as aforementioned, overshooting and ringing (clearly perceived in the 1D profile). These effects are much milder in the L1 reconstruction.

\begin{figure*}
    \centering
    \subfloat[Depth frames $i-1$, $i$, and $i+1$]{
        \begin{tabular}{ccc}
            \includegraphics[width=\mylength]{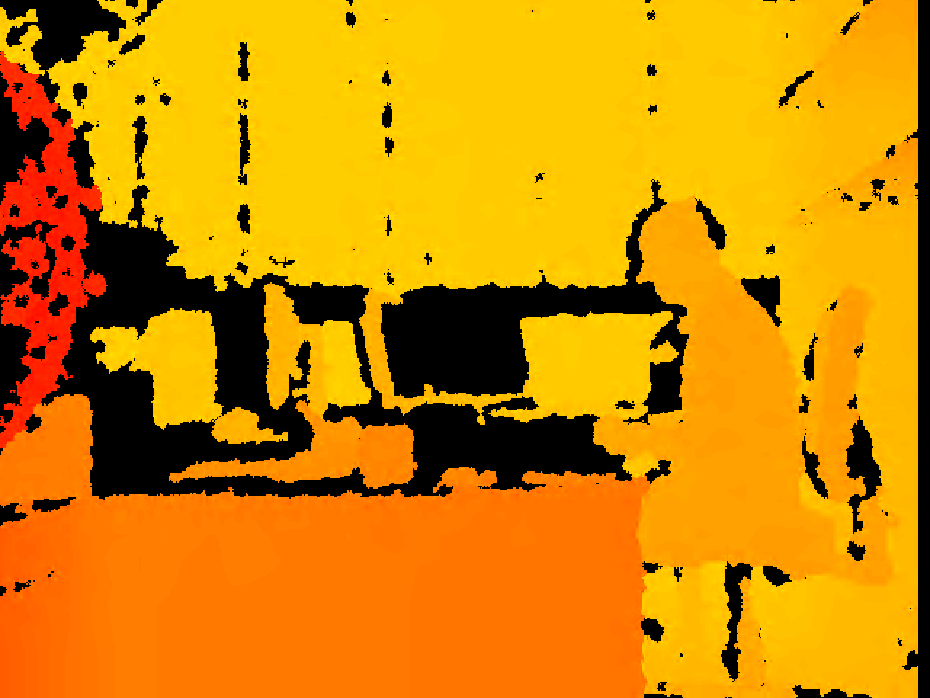}
            \includegraphics[width=\mylength]{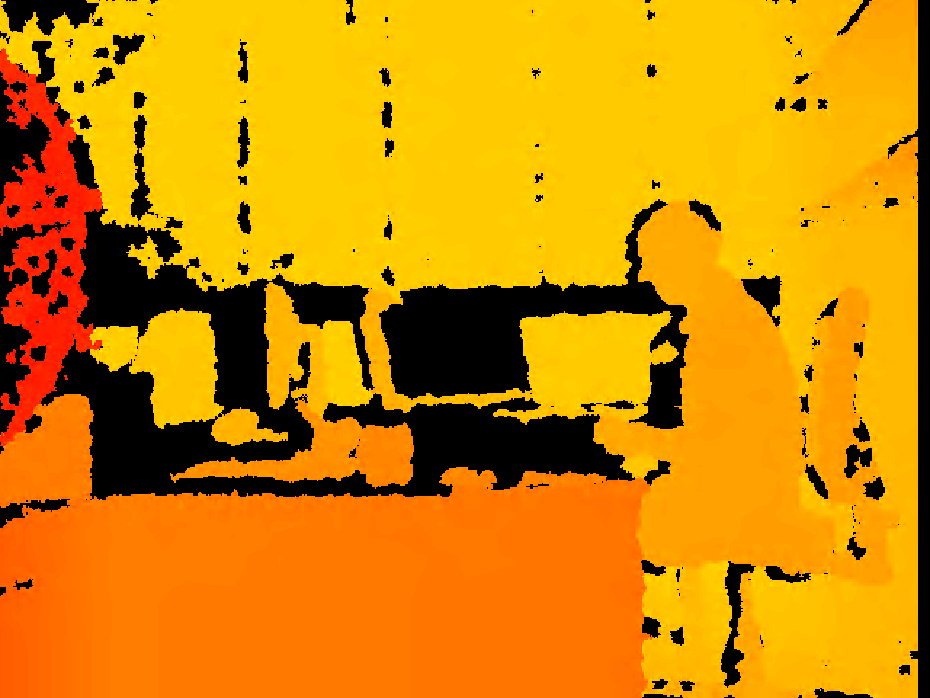}
            \includegraphics[width=\mylength]{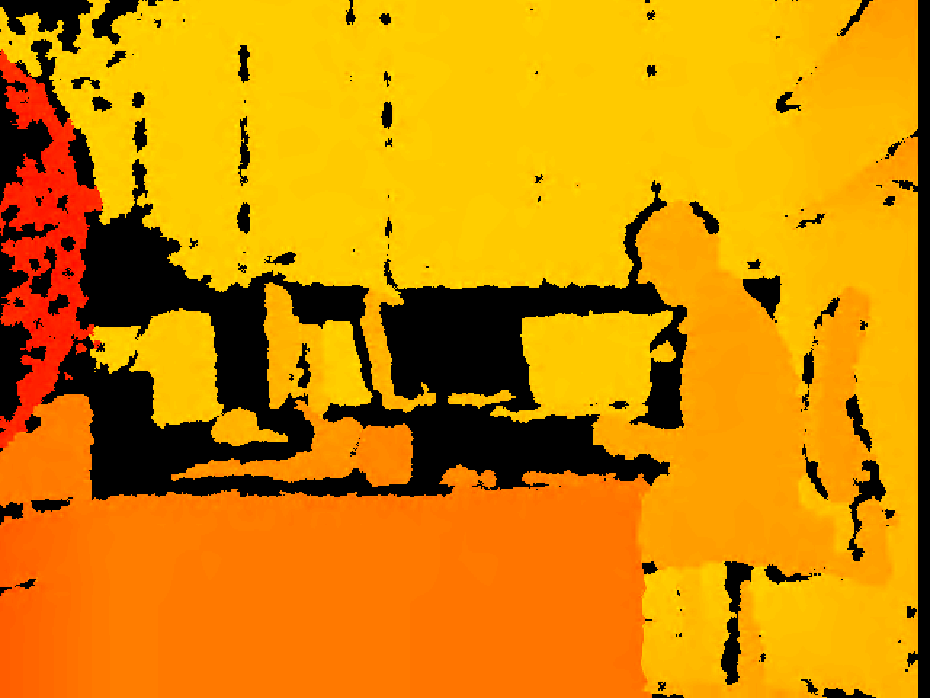}
        \end{tabular}
        \label{fig:kinect3D_original}
    }
    
    \centerline{
        \hfill
        \subfloat[Robust L2 spline]{
            \includegraphics[width=\mylength]{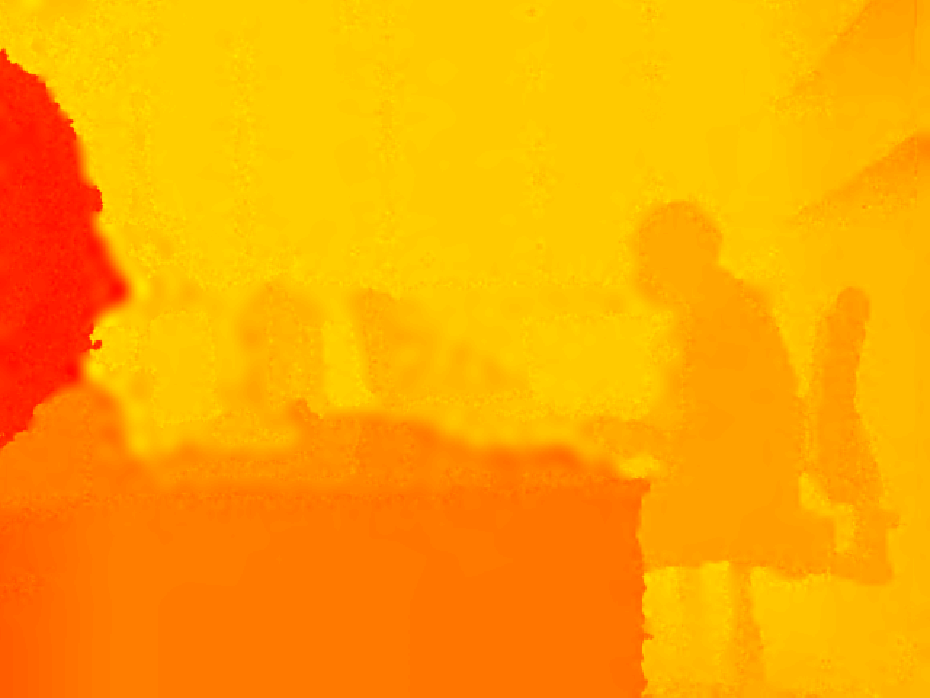}
        }
        \hfill
        \subfloat[L1 spline]{
            \includegraphics[width=\mylength]{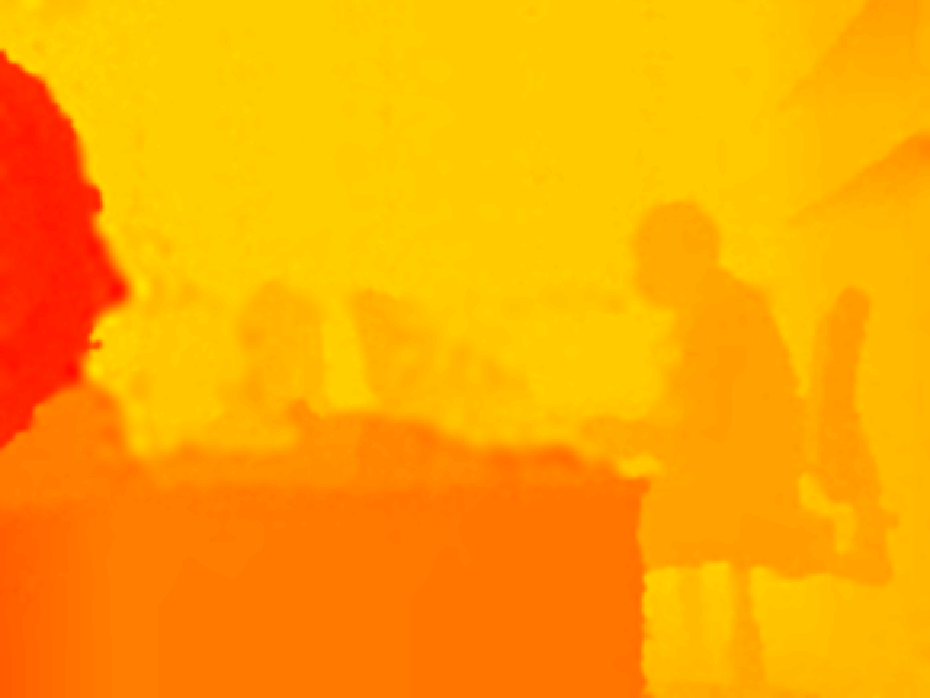}
        }
        \hfill
    }
    
    \subfloat[Profile of a single row. Left, original values; right, original and reconstructed values.]{
        \centerline{
            \hfill
            \includegraphics[width=0.45\textwidth]{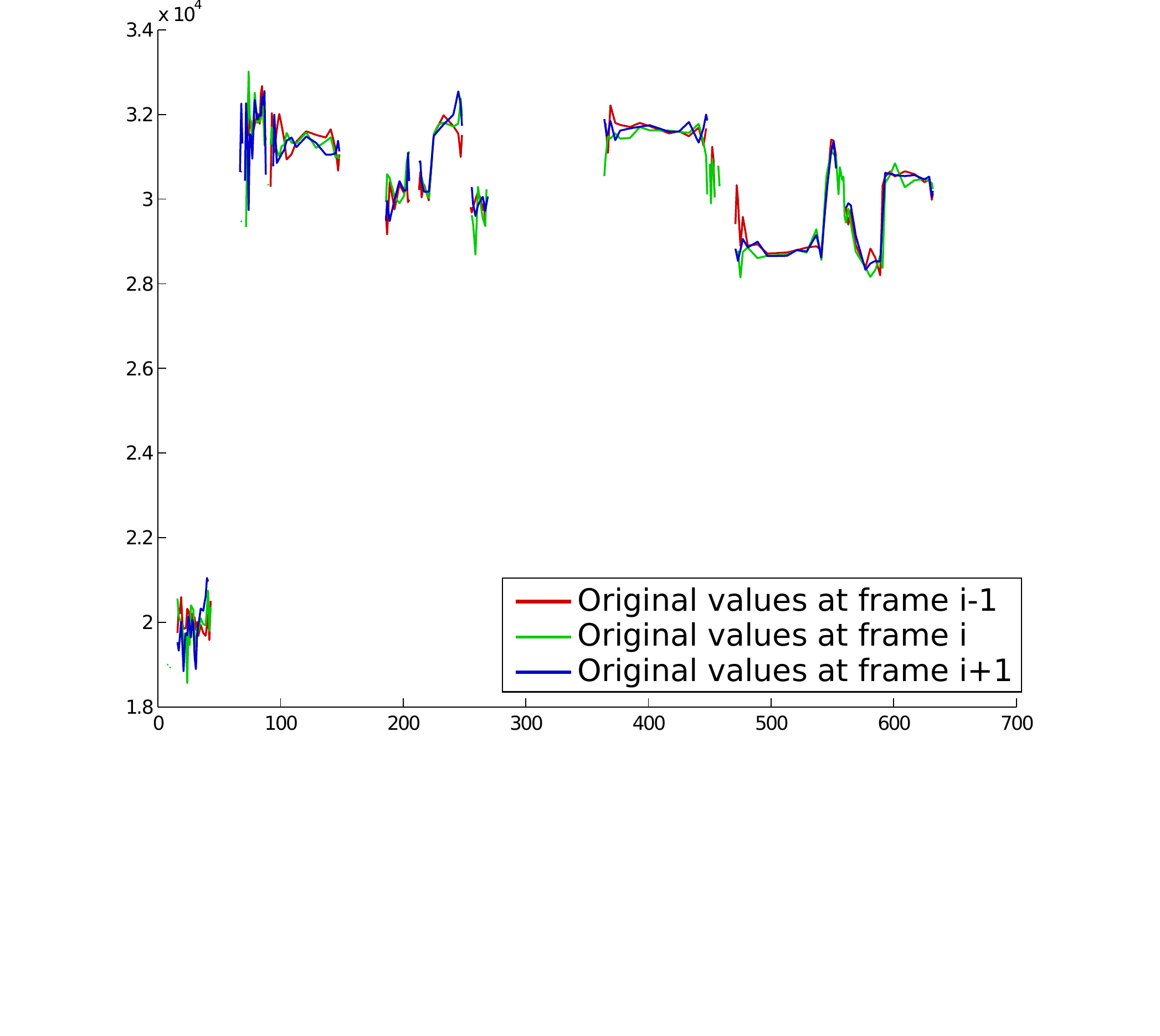}
            \hfill
            \includegraphics[width=0.45\textwidth]{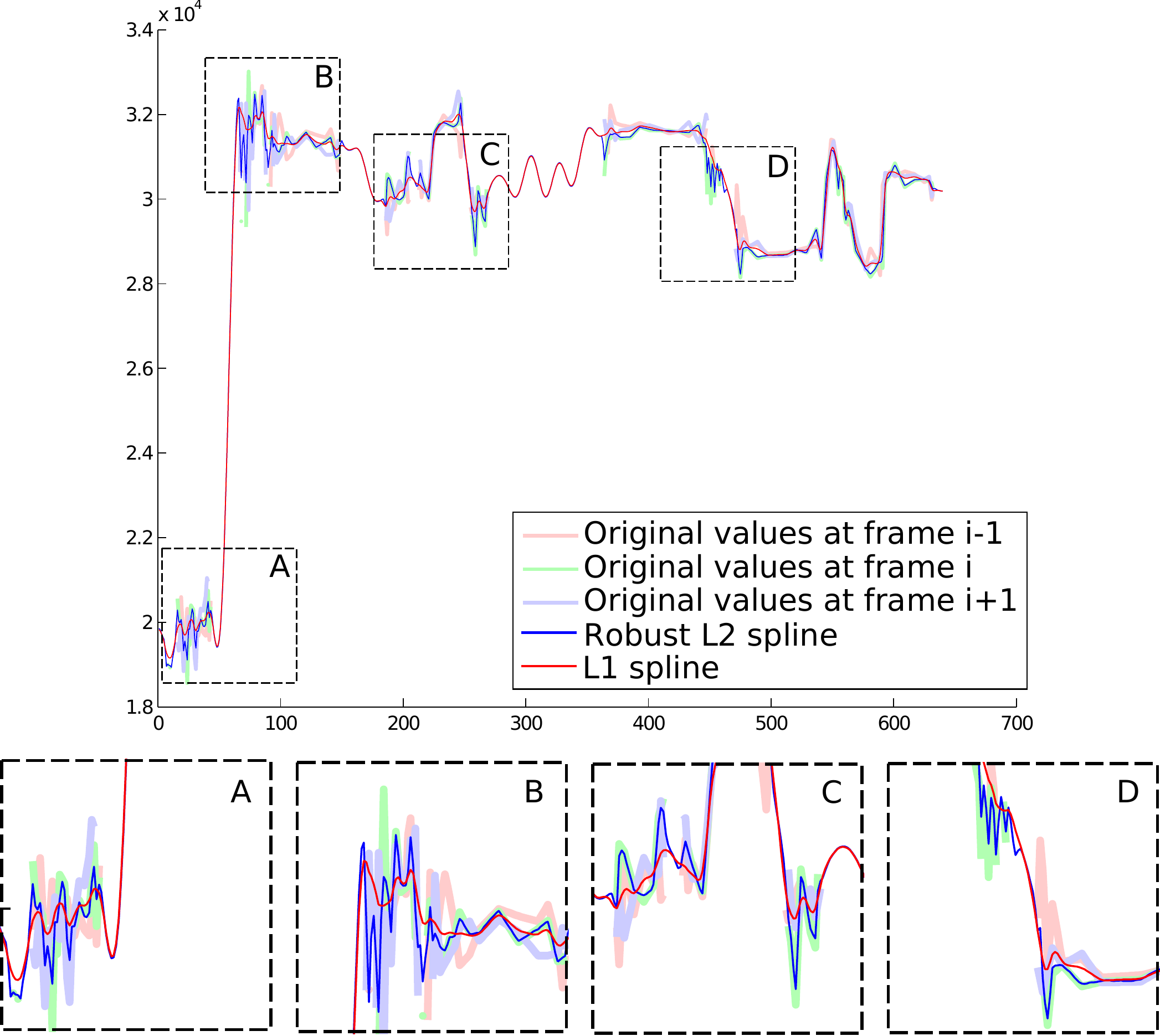}
            \hfill
        }
    }

    \caption{Smoothing and interpolating using multiple frames. To process frame $i$, we build a $640 \times 480 \times 3$ tensor using frames $i-1$, $i$, and $i+1$, which we then smooth/interpolate with a 3D robust L2 or L1 spline. Since the data presents several ``jumps,'' the robust L2 spline must under-smooth the data to be able to fit it correctly. The L1 spline presents a good trade-off between fitting and smoothing.}
    \label{fig:kinect3D}
\end{figure*}

In Fig.~\ref{fig:kinect3D_original} we can clearly observe that the position of the missing data is not consistent across frames. We can integrate data from several frames to achieve more accurate interpolations, by performing 3D reconstructions.
Thus, in this example, we treat depth data as a 3D signal (2D + time), by considering three consecutive frames. The data dimensionality is then $640 \times 480 \times 3$. A full depth video can be smoothed by using 3D splines as a sliding-window type of filter. The robust L2 spline again presents a noisier behavior and with significative overshooting. On the other hand, the L1 spline is much smoother in smooth areas while correctly preserving abrupt transitions.

\noindent\textbf{Running times.}
We present in Table~\ref{tab:runningTime} the running-time and number of iterations until convergence for every example in this work. The time of the robust L2 and the L1 splines is comparable. All code is written in pure Matlab, with no C++ or mex optimizations. All experiments were run on  a MacBook Pro with a 2.7GHz Intel Core i7 processor.
Finally, note that in most cases the algorithm converges in less than twenty outer-iterations (recall that $\varepsilon=10^{-3}$). In the example in Fig.~\ref{fig:spline_2D_n30}, the maximum number of iterations (100) is reached with a final error of $10^{-2.8}$.

\begin{table}
    \caption{Execution times (in seconds) and number of iterations until convergence of the proposed algorithm for the different experiments performed in this work.}
    \label{tab:runningTime}
    \centering
    \begin{tabular}{llcccc}
        \toprule
        && \multirow{2}{*}{Size} &
        Robust L2 spline & \multicolumn{2}{c}{L1 spline} \\
        \cmidrule{5-6}
        &&& Time & Time & Iters. \\
        \midrule
        \multirow{7}{*}{\begin{minipage}{0.45in}Complete\\data\end{minipage}}
        & Fig.~\ref{fig:spline_1D_16} & $2^{20}$ & 4.931 & 3.590 & 7 \\
        & Fig.~\ref{fig:spline_2D_n30} & $256 \times 256$ & 0.541 & 0.982 & 100 \\
        & Fig.~\ref{fig:temp} & 163 & 0.007 & 0.045 & 72 \\
        & Fig.~\ref{fig:load} & 501 & 0.010 & 0.010 & 11 \\
        & Fig.~\ref{fig:square} & $2^{10}$ & 0.035 & 0.007 & 6 \\
        & Fig.~\ref{fig:119082_original} & $321 \times 481$ & 1.167 & 0.758 & 17 \\
        & Fig.~\ref{fig:119082_noise} & $321 \times 481$ & 1.143 & 0.873 & 20 \\
        \midrule
        \multirow{2}{*}{\begin{minipage}{0.4in}Missing\\data\end{minipage}}
        & Fig.~\ref{fig:kinect2D} & $480 \times 640$ & 0.911 & 0.266 & 2 \\
        & Fig.~\ref{fig:kinect3D} & $480 \times 640 \times 3$ & 7.511 & 2.202 & 3\\
        \bottomrule
    \end{tabular}
\end{table}

\section{Conclusions}
\label{sec:conclusions}

We have presented a new method for robustly smoothing regularly sampled data. We do this with modified splines, where we replace the classical L2-norm in the fitting term by an L1-norm. This automatically handles outliers, thus obtaining a robust approximation.

We also presented a new technique, using split-Bregman iteration, for solving the resulting optimization problem. The algorithm is extremely simple and easy to code. The method converges very quickly and has a small memory footprint. It also makes extensive use of the DCT, thus being straightforward to implement in GPU. These characteristics make this method very suitable for large-scale problems.

\section*{Acknowledgment}

Work partially supported by NSF, ONR, NGA, ARO, DARPA, and NSSEFF.
We thank Dr.~Gonzalo Mateos for kindly providing the power consumption dataset.

\bibliographystyle{plain}
\bibliography{mtepper}

\end{document}